\documentclass{amsart}
\usepackage{graphicx}
\usepackage{amsmath}
\usepackage{color}

\newtheorem{thm}[equation]{Theorem}
\newtheorem{pro}[equation]{Proposition}
\newtheorem{cor}[equation]{Corollary}
\newtheorem{lem}[equation]{Lemma}

\theoremstyle{definition}
\newtheorem{exa}[equation]{Example}

\newtheorem{DEF}[equation]{Definition}
\newtheorem{rem}[equation]{Remark}

\usepackage{amsfonts,amssymb,amscd,amsmath,enumerate,verbatim,calc}

\def\aa{\mathcal A}

\def\bbbf{\mathbb{F}}

\def\LL{\mathcal{L}}

\def\lam{\lambda}
\def\Lam{\Lambda}

\def\u{{\mathcal U}}
\def\v{{\mathcal V}}

\def\l{L}

\def\g{G}

\def\vv{{\mathcal V}}
\def\uu{{\mathcal U}}
\def\aa{{\mathcal A}}

\def\lam{\Lambda_I^{\bar{i}}}

\def\lam{\lambda}
\def\Lam{\Lambda}
\def\ad{\hbox{ad}}
\begin{document}

\markboth{On the structure of graded $3$-Lie-Rinehart algebras} {V. Khalili}

\date{}

\centerline{\bf On the structure of graded $3$-Lie-Rinehart algebras }

\vspace{.5cm}\centerline{Valiollah Khalili\footnote[1]{Department
of mathematics, Faculty of sciences, Arak University, Arak 385156-8-8349, Po.Box: 879, Iran. 
 V-Khalili@araku.ac.ir\\
\hphantom{ppp}2000 Mathematics Subject Classification(s):17B05, 17B22, 17B60, 17A60.
\\
\hphantom{ppp} Keywords:$3$-Lie-Rinehart algebra, graded algebra and structure theory. }\;\;}

\vspace{1cm} \noindent ABSTRACT: We study the structure of a  graded $3$-Lie-Rinehart algebra $\LL$  over an associative  and commutative graded algebra $A.$ For $G$ an abelian group,   we show that if $(\LL, A)$  is a tight $G$-graded $3$-Lie-Rinehart  algebra, then $\LL$ and $A$ decompose as  $\LL =\bigoplus_{i\in I}\LL_i$ and $A =\bigoplus_{j\in J}A_j,$ where any $\LL_i$ is a non-zero  graded ideal of $\LL$ satisfying $[\LL_{i_1}, \LL_{i_2}, \LL_{i_3}]=0$ for any $i_1, i_2, i_3\in I$  different from each other,  and any $A_j$ is a non-zero  graded  ideal of $A$  satisfying $A_j A_l=0$ for any $l, j\in J$ such that $j\neq l,$ and both decompositions satisfy that for any $i\in I$  there exists a unique $j\in J$ such that  $A_j \LL_i\neq 0. $ Furthermore, any $(\LL_i, A_j)$  is a graded $3$-Lie-Rinehart algebra. Also,  under certain conditions, it is shown that the above decompositions of $\LL$ and $A$ are by means of the family of their, respectively, graded simple ideals.
 
\vspace{1cm} \setcounter{section}{0}
\section{Introduction}\label{introduction}

The notion of Lie-Rinehart algebra plays an important role in many branches of
mathematics. They are algebraic analogs of Lie algebroids.  The idea of this notion first introduced by Herz \cite{H} as pseudo-Lie algebras, then studied by Palais \cite{P} under the name "d-Lie ring". Lie-Rinehart structures have been the subject of extensive studies, such as in relation to  differential geometry \cite{Ri},  differential Galois theory \cite{H3}, symplectic geometry \cite{LSB, LSBC}, Poisson structures \cite{Pi}, various kinds of quantizations \cite{H0, H00}, Lie groupoids and Lie algebroids \cite{M, S, SZ}. For a very extensive survey of those topics, the reader can be found in \cite{Ch, D, H1, H2, K, SC}.

The study of gradings on Lie algebras begins in the 1933 by Jordan’s
work \cite{J}, with the purpose of formalizing Quantum Mechanics. Since then,
many papers describing different physic models by means of graded Lie type
structures have appeared, being remarkable the interest on these objects in
the last years. It is worth mentioning that the so-called techniques of connection of roots had long been introduced by Calderon, Antonio J, on split Lie algebras with symmetric root systems in \cite{C1}. For instance, in  reference \cite{C2} the author studied the structure of arbitrary graded Lie algebras, being extended to the framework of graded Lie superalgebras in \cite{CS1} by the technique of connections of elements in the support of the graing. Recently, in \cite{C3, C4, Kh}, the structure of arbitrary  graded commutative algebras, graded Lie triple systems  and graded $3$-Leibniz algebras  have been determined by the  connections of the support of the grading. 

A $3$-Lie-Rinehart algebra is a triple $(\LL, A, \rho),$ where $\LL$ is a $3$-Lie  algebra, $A$ is a commutative, associative  algebra, $\LL$ is an $A-$module, $(A, \rho)$ is a  $\LL-$module in such a way that both structures are related in an appropriate way. Our goal in this work is to study the inner structure of arbitrary graded $3$-Lie-Rinehart  algebras by the  developing technique of connections of elements in the supports of the graing of $\LL$ and  $A.$  The finding of the present paper is an improvement and extension of the work on graded  Lie-Rinehart algebras in \cite{BCNS}.

The article is organized as follows; In Section 2, we recall the definition of  $3$-Lie-Rinehart algebras  and introduced a class of graded $3$-Lie-Rinehart algebra by means of the abelian group $G.$ In Section 3, as a second step, we extend the techniques of connections in the support of the graing for graded Lie algebras in \cite{C2} to the framework of graded $3$-Lie-Rinehart  algebra $(\LL, A).$  In Section 4,  we get, as a third  step, a decomposition of $A$ as direct sum of adequate ideals. We also charactrized the relation between the decomposition of $\LL$ which obtained in section 3 and the given decomposition of $A.$  Section 5 is devoted to show that, under mild conditions, the given decompositions of $\LL$ and $A$ are by means of the family of their, corresponding, graded simple ideals.
 
Throughout this paper, algebras and vector spaces are over a field $\bbbf$ of characterestic zero, and $A$ denotes an associative and commutative algebra over $\bbbf.$ We also
consider an  abelian group $G$ with unit element 1.

\section{Preliminaries} \setcounter{equation}{0}\

In this section, we recall definitions and some results on $3$-Lie-Rinehart algebras and alsointroduced a class of graded $3$-Lie-Rinehart algebra by means of the abelian group $G.$

\begin{DEF}\label{3L}\cite{F} A $3$-Lie algebra consists of a  vector space $\LL$
together with a trilinear map $[., ., .] : \LL\times\LL\times\LL\longrightarrow\LL$ such that the following
conditions are satisfied:
\begin{itemize}
\item[(i)] skew symmetry:
$[x_1, x_2, x_3]= -[x_2, x_1, x_3]= -[x_1, x_3, x_2];$
\item[(ii)] fundamental identity:
\begin{eqnarray}\label{0}
[[x_1, x_2,  x_3], y_1, y_2]&=& [[x_1, y_1, y_2], x_2, x_3]+[[x_2, y_1, y_2], x_3,  x_1]\\
\nonumber&+&[[x_3, y_1, y_2], x_1, x_2],
\end{eqnarray}
\end{itemize}
for all  elements $x_1, x_2,  x_3, y_1, y_2\in\LL.$
\end{DEF}

\begin{DEF}\cite{Ka} Let $(\LL, [., ., .])$ be a $3$-Lie   algebra,  $V$ be a  vector space and $\rho:\LL\times\LL \longrightarrow gl(V)$ be a linear mapping.
Then $(V,\rho)$ is called a representation of $\LL$ or $V$ is an $\LL$-module  if the following two conditions hold:
\begin{itemize}
\item[(i)] $[\rho(x_1, x_2), \rho(x_3, x_4)] = \rho([x_1, x_2, x_3], x_4) -\rho([x_1, x_2, x_4], x_3),$
\item[(ii)] $\rho([x_1, x_2, x_3], x_4)= \rho(x_1, x_2)\rho(x_3, x_4) +\rho(x_2, x_3)\rho(x_1, x_4)
+\rho(x_3, x_1)\rho(x_2, x_4),$
\end{itemize}
for all  elements $x_1, x_2,  x_3, x_4\in\LL.$
\end{DEF}

Next, define
$$
\ad : \LL\times\LL\longrightarrow gl(\LL);~~~ad(x, y)z=[x,y,z] ,~~\forall x, y, z\in\LL.
$$
Tanks to fundamental identity, $(\LL, \ad)$ is a representation of the  $3$-Lie   algebra $\LL,$ and it is called the adjoint representation of $\LL.$ One can see that $\ad(\LL, \LL)$ is a Lie algebra which is called inner derivation of $\LL.$ We also have by fundamental identity,
$$
[\ad(x_1, x_2), ad(y_1, y_2)]=\ad([x_1, y_1, x_2], y_2)+\ad(x_2, [x_1, y_1, y_2]).
$$

\begin{DEF}\cite{BL}  Let   $(\LL, [., ., .])$ be a $3$-Lie   algebra, $\LL$ be  an $A$-module and $(A, \rho)$ be an $\LL$-module. If $\rho(\LL, \LL)\subset Der(A)$ and,
\begin{eqnarray}\label{4}
[x, y, a z] = a[x, y, z]+\rho(x, y) a z,~~~\forall x, y, z\in\LL,~~\forall a\in A,
\end{eqnarray}
\begin{eqnarray}\label{5}
\rho(ax, y) = \rho(x, ay)= a\rho(x, y),~~~\forall x, y\in\LL,~~\forall a\in A,
\end{eqnarray}	
then $(\LL, A, [., ., .], \rho)$	is called a {\em $3$-Lie-Rinehart  algebra}.
\end{DEF}

\begin{rem}
\begin{itemize}
\item[(1)]  If $\rho=0$ then $(\LL, A, [., ., .])$ is a $3$-Lie $A$-algebra.
\item[(2)] Let $(\g, [., .])$ be a  Lie algebra and an $A-$module. Let $(A, \rho)$ be a $\g$-module. If $\rho(\g)\subset Der(A)$  and
$$
[x, ay]=)a[x, y]+\rho(x)ay,~~~~\rho(ax)=a\rho(x),~~~~~\forall x, y\in\g,  ~  \forall a\in A,
$$
then $(\g, [., .], A, \rho)$ a Lie-Rinehart  algebra.
\end{itemize} 
\end{rem}

\begin{exa} We recall that given a Lie  algebra analogues of  trace, 
	we can construct a $3$-Lie  algebra. Let $(\LL, [., . ])$ be a Lie  algebra and
	$\tau : \LL\longrightarrow\bbbf$ an  linear form. We say that $\tau$ is a  trace of $\LL$ if $\tau([., .]) = 0.$ For
	any $x_1, x_2, x_3\in\LL,$ we define the $3$-ary bracket by
	\begin{equation}\label{500}
	[x_1, x_2, x_3]_\tau = \tau(x_1)[x_2, x_3] - \tau(x_2)[x_1, x_3]+\tau(x_3)[x_1, x_2]. 
	\end{equation}	
	Then $(\LL, [., ., .]_\tau)$ is a $3$-Lie  algebra (see  \cite{A} for super algebras).
	Next, we begin by constructing
	$3$-Lie-Rinehart  algebras starting with a Lie-Rinehart  algebras. Let $(\LL, A, [., .], \rho)$ be a Lie-Rinehart algebra and $\tau$ is a  trace of $\LL.$ If the condition
	$$
	\tau(ax)y = \tau(x)ay,
	$$ 
	is satisfied for any  $x, y\in\LL, a\in A$, then $(\LL, A, [., ., .]_\tau, \rho_\tau)$ is a $3$-Lie-Rinehart  algebra, where $[., ., .]_\tau$ is defined as Eq.(\ref{500}) and $\rho_\tau$ is defined by
	$$
	\rho_\tau :\LL\times\LL\longrightarrow gl(\LL) ;~~\rho_\tau(x, y)=\tau(x)\rho(y)-\tau(y)\rho(x),~~\forall x, y\in\LL.
	$$
	(see Theorem 2.1 in \cite{BhCMS} for super algebras)	
\end{exa}

\begin{DEF}\label{ideal} Let $(\LL, A, [., ., .], \rho)$	be a $3$-Lie-Rinehart  algebra.
\begin{itemize}
\item[(1)] If $S$ is a $3$-Lie  subalgebra  of  $\LL$ satisfying $A S\subset S$ and 
$$
(S, A, [., ., .]|_{S\times S}, \rho|_{S\times S})
$$
(which is a  $3$-Lie-Rinehart  algebra)  is called a {\em subalgebra} of the   $3$-Lie-Rinehart  algebra  $(\LL, A, [., ., .], \rho).$

\item[(2)]  If $I$ is a $3$-Lie  ideal  of  $\LL$ satisfying  $\rho(I, I)(A)(\LL)\subset  I$ and 
$$
(I, A, [., ., .]|_{I\times I}, \rho|_{I\times I})
$$
(which  is a  $3$-Lie-Rinehart  algebra) is called an {\em ideal} of the   $3$-Lie-Rinehart  algebra  $(\LL, A, [., ., .], \rho).$ 
\item[(3)]  We also  say that  $(\LL, A, [., ., .], \rho)$ is {\em simple} if $[\LL, \LL, \LL]\neq 0,~AA\neq 0,  A\LL \neq 0$ and its only ideals are $\{0\}$, $\LL$ and $\ker \rho.$ 
\end{itemize}
\end{DEF}

For a $3$-Lie-Rinehart  algebra  $(\LL, A, [., ., .], \rho),$ we denote
$$
Ann(A) :=\{a\in A : aA=0\},~~~~\hbox{and~~}~~~~Ann_{\LL}(A) :=\{a\in A : ax=0, \forall x\in\LL\},\\
$$
the annihilator of the commutative and
associative algebra $A$ and the annihilator of $A$ in $\LL.$ We also denote
$$
Z_\rho(\LL) :=\{x\in\LL : [x, \LL, \LL]=0,~\rho(x, \LL)=0\},~~~~\hbox{and~~}~~~~\ker \rho :=\{x\in\LL : \rho(x, \LL)=0\},\\
$$
the center of  $3$-Lie-Rinehart  algebra $(\LL, A, [., ., .], \rho)$ and the kernel of the representation.  Note that $Z_\rho(\LL)=\ker \rho\cap Z(\LL),$  and by Theorem 2.3 in \cite{BL}, $Ann_{\LL}(A)$ is an ideal of $A$ and $Z_\rho(\LL)$ is an ideal of  $3$-Lie-Rinehart  algebra $(\LL, A, [., ., .], \rho).$

\begin{DEF}\label{g3l} Let $\LL$ be a $3$-Lie algebra. It is said that $\LL$ is graded by means
of an abelian group $G$ if it decomposes as the direct sum of linear subspaces
$$
\LL=\bigoplus_{g\in G} \LL_g,
$$
where the homogeneous components satisfy $[\LL_g, \LL_h, \LL_k] \subset \LL_{ghk}$ for any $g, h, k \in	G$ (denoting by juxtaposition the product and unit element 1 in $G$). 
\end{DEF}

Note that  split $3$-Lie algebras \cite{CP} and graded Lie triple systems \cite{C4} are examples of graded $3$-Lie algebras.

\begin{DEF}\label{maindef} We say that a  $3$-Lie-Rinehart  algebra  $(\LL, A)$  is a graded algebra, by means
of the abelian group $G,$ if $\LL$ is a $G$-graded $3$-Lie algebra as in Definion \ref{g3l} and the algebra $A$ is a $G$graded (commutative and associative) algebra in the sense that
$A$ decomposes as $A=\bigoplus_{h\in G}A_h,$ with $A_g A_h\subset A_{gh},$ satisfying
\begin{equation}\label{00}
A_h\LL_g\subset\LL_{hg}
\end{equation}
\begin{equation}\label{000}
\rho(\LL_g, \LL_{g'})(A_h)\subset A_{gg'h},
\end{equation}
for any $g, g',  h\in G.$
\end{DEF}
 split  $3$-Lie-Rinehart  algebra is example of graded  $3$-Lie-Rinehart  algebra. So this paper extends the results obtained in \cite{Kh1}.

As it is usual in the theory of graded algebras, the  regularity conditions will be understood in the graded sense compatible with the $3$-Lie-Rinehart algebra structure. That
is, a  $3$-Lie-Rinehart graded subalgebra (or graded ideal) of $(\LL, A)$ is a graded linear subspace $S$ (or $I$) as in Definition \ref{ideal}.   More precisely, $S$ (or $I$) splits as $S=\bigoplus_{g\in G} S_g,~~S_g=S\cap \LL_g,$ similarly for $I.$ Also we will say that $(\LL, A)$ is a gr-simple  $3$-Lie-Rinehart algebra if $[\LL, \LL, \LL]\neq 0$ and its only graded ideals are $\{0\},  \LL$ and $\ker \rho.$  

 We denote the $G$-support of the grading in $\LL$ and in $A$ to the sets
$$ 
\Sigma^1=\{g\in G\setminus\{1\}~:~\L_g\neq 0\}~~~\hbox{and }~~~\Lambda^1=\{h\in G\setminus\{1\}~:~A_h\neq 0\},
$$
respectively. If  $(\LL, A)$  is a graded $3$-Lie-Rinehart algebra then we can rewrite
$$
\LL=\LL_1\oplus(\bigoplus_{g\in \Sigma^1}\LL_g)~~~\hbox{and }~~~A=A_1\oplus(\bigoplus_{h\in \Lambda^1}A_h).
$$

\section{Connections in $\Sigma^1$ and decompositions} \setcounter{equation}{0}\
In this section, we begin by developing the techniques
of connections  in $\Sigma^1.$  Let  $(\LL, A)$ be a graded  $3$-Lie-Rinehart  algebra,  with the decomposition
$$
\LL=\LL_1\oplus(\bigoplus_{g\in \Sigma^1}\LL_g)~~~\hbox{and }~~~A=A_1\oplus(\bigoplus_{\lam\in \Lambda^1}A_{\lam}),
$$
 and  with  the $G$-supports $\Sigma^1$ and $\Lambda^1,$ respectively. 

We define
$$ 
\Sigma^{-1}=\{g^{-1}~:~g\in \Sigma^1\},~~~\hbox{and }~~~\Lambda^{-1}=\{\lam^{-1}~:~\lam\in\Lambda^1\}.
$$
Let us denote
$$
\Sigma=\Sigma^1\cup\Sigma^{-1},~~~\hbox{and }~~~\Lam = \Lambda^1\cup\Lambda^{-1}.
$$
\begin{DEF}\label{conn1}
Let $g, h$ be two elements in $\Sigma^1.$ We say that {\em $g$
is $\Sigma^1$-connected to $h$} if there exists a
family $\{g_1, g_2, g_3, ..., g_{2n+1}\}\subset\Sigma\cup\Lam\cup\{1\},$
satisfying the following conditions;\\

\begin{itemize}
\item[(1)] $g=g_1,$
\item[(2)] $ \{ g_1g_2g_3, g_1g_2g_3g_4g_5, ..., g_1g_2g_3... g_{2n-1}\}\subset \Sigma,$
\item[(3)] $g_1g_2g_3 ... g_{2n+1}\in\{h, h^{-1}\}.$
\end{itemize}
The family $\{g_1, g_2, g_3, ..., g_{2n+1}\}$ is called a
{\em $\Sigma^1$-connection} from $g$ to $h.$
\end{DEF}

The next result shows that the $\Sigma^1$-connection relation is an equivalence relation. Its
proof is analogous to the one for graded Lie triple system   in \cite{BCNS}, Proposition 3.1.

\begin{pro}\label{rela}
The relation $\sim_{\Sigma^1}$ in $\Sigma^1$ defined by
$$
g\sim_{\Sigma^1} h~\hbox{ if and only if}~ g~\hbox{ is $\Sigma^1$-connected to}~ h,
$$
is an equivalence relation.
\end{pro}

By the  Proposition \ref{rela}, we can consider the equivalence relation
in $\Sigma^1$ by the connection relation $\sim_{\Sigma^1}.$  So we denote
by
$$
\Sigma^1/\sim_{\Sigma^1} :=\{[g] : g\in\Sigma^1\},
$$
where $[g]$ denotes  the set of elements of $\Sigma^1$  which are
connected to $g.$

 Clearly, if $h\in [g],$ then $h^{-1}\in [g].$ 

\begin{rem}\label{31} For $g', g''\in\Sigma\cup\Lam\cup\{1\},$ if $h\in [g]$ and $gg'g''\in\Sigma^1$ then $h\sim_{\Sigma^1}gg'g''.$ Indeed, the family $\{g, g', g''\}$ is a $\Sigma^1$-connection from $g$ to $gg'g''.$ Now, taking into account $g\sim_{\Sigma^1} h$ and Proposition \ref{rela},  we get $h\sim_{\Sigma^1}gg'g''.$
\end{rem}

Our next goal  is to associate an adequate ideal
$I_{[g]}$ of $\LL$ to any $[g].$ For a fixed  $g\in\Sigma^1,$ we define
 \begin{eqnarray}\label{38}
\LL_{1, [g]} :=\big(\sum_{h\in[g]\cap\Lam^1} A_{h^{-1}}\LL_h\big)+\big(\sum_{h, k\in[g]}[\LL_h, \LL_k, \LL_{(hk)^{-1}}]\big)\subset\LL_1.
\end{eqnarray}
Next, we define
\begin{equation}\label{38'}
\vv_{[g]}
:=\bigoplus_{h\in[g]}\LL_h.
\end{equation}
Finally, we denote by $I_{[g]}$ the direct sum of the two 
subspaces above, that is,
\begin{equation}\label{38''}
I_{[g]} :=\LL_{1, [g]}\oplus\vv_{[g]}.
\end{equation}

\begin{pro}\label{subalg} For any $[g]\in\Sigma^1/\sim_{\Sigma^1},$
 the following assertions hold.
 \begin{itemize}
 \item[(1)] $[I_{[g]}, I_{[g]}, I_{[g]}]\subset I_{[g]}.$
 \item[(2)]$A I_{[g]}\subset I_{[g]} $
 \end{itemize}
\end{pro}
\noindent {\bf Proof.} (1) By Eq. (\ref{38''})  we have
\begin{eqnarray}\label{200}
\nonumber[I_{[g]}, I_{[g]},  I_{[g]}]&=&[\LL_{1, [g]}\oplus\vv_{[g]}, \LL_{1, [g]}\oplus\vv_{[g]}, 
\LL_{1, [g]}\oplus\vv_{[g]}]\\
&\subset&[\LL_{1, [g]}, \LL_{1, [g]},   \LL_{1, [g]}]+[\LL_{1, [g]}, \LL_{1, [g]}, \vv_{[g]}]+[\LL_{1, [g]}, \vv_{[g]}, \LL_{1, [g]}]\\
\nonumber&+&[\LL_{1, [g]}, \vv_{[g]}, \vv_{[g]}]+[\vv_{[g]}, \LL_{1, [g]},   \LL_{1, [g]}]+[\vv_{[g]}, \LL_{1, [g]}, \vv_{[g]}]\\
\nonumber&+&[\vv_{[g]}, \vv_{[g]},  \LL_{1, [g]}]+[\vv_{[g]}, \vv_{[g]}, \vv_{[g]}].
\end{eqnarray}
Since $ \LL_{1, [g]}\subset\LL_1$ and by the skew symmetry of  trilinear map, we clearly have
\begin{equation}\label{38.5}
[\LL_{1, [g]}, \LL_{1, [g]}, \vv_{[g]}]+[\LL_{1, [g]}, \vv_{[g]}, \LL_{1, [g]}]+[\vv_{[g]}, \LL_{1, [g]},   \LL_{1, [g]}]\subset \vv_{[g]}.
\end{equation}

Consider now the   summand $[\LL_{1, [g]}, \vv_{[g]}, \vv_{[g]}]$ in (\ref{200}). By $\LL_{1, [g]}\subset\LL_1,$ we have  
\begin{equation}\label{39}
[\LL_{1, [g]}, \vv_{[g]}, \vv_{[g]}]\subset [\LL_1, \vv_{[g]}, \vv_{[g]}].
\end{equation}
 Suppose there exist $g', g''\in[g]$ such that $[\LL_1, \LL_{g'}, \LL_{g''}]\neq 0.$ Then $g'\in\Sigma^0$ and $g'g''\in\Sigma^1\cup\{1\}.$  If $g''=g'^{-1},$ clear that $[\LL_1, \LL_{g'}, \LL_{g''}]=[\LL_1, \LL_{g'}, \LL_{g'^{-1}}]\subset \LL_{1, [g]}.$ Otherwise, if $g''\neq g'^{-1},$  and  $\{g_1, g_2, g_3, ..., g_{2n+1}\}$ is  a
$\Sigma^1$-connection from $g$ to $h.$ Then $\{g_1, g_2, g_3, ..., g_{2n+1}, 1, k\}$  is  a
$\Sigma^1$-connection from $g$ to $hk$ in case $g_1 g_2 g_3...g_{2n+1}=h$  and  $\{g_1, g_2, g_3, ..., g_{2n+1}, 1, k^{-1}\}$   in case $g_1 g_2 g_3...g_{2n+1}=h^{-1}.$  Hence, $hk\in[g].$ Taking into account Eq. (\ref{39}), we get $[\LL_{1, [g]}, \vv_{[g]}, \vv_{[g]}]\subset \vv_{[g]}.$ By  the skew symmetry of  trilinear map, we get
\begin{equation}\label{41}
[\LL_{1, [g]}, \vv_{[g]}, \vv_{[g]}]+[\vv_{[g]}, \LL_{1, [g]}, \vv_{[g]}]+[\vv_{[g]}, \vv_{[g]},  \LL_{1, [g]}]\subset \vv_{[g]}.
\end{equation}

Consider now the   summand $[\vv_{[g]}, \vv_{[g]}, \vv_{[g]}]$ in (\ref{200}).  Suppose there exist $h, k, l\in[g]$ such that $[\LL_h, \LL_k, \LL_l]\neq 0.$ Then $hk\in\Sigma^0\cup\{1\}$ and $hkl\in\Sigma^1\cup\{1\}.$ If either $h=k^{-1}$ or $hkl=1,$ then 
$$
[\LL_h, \LL_k, \LL_l]=\LL_l\subset\vv_{[g]}~~~~~\hbox{or}~~~~~[\LL_h, \LL_k, \LL_l]\subset\LL_{1, [g]}.
$$
 Otherwise, if $hk\in\Sigma^0$  and  $hkl\in\Sigma^1,$ and $\{g_1, g_2, g_3, ..., g_{2n+1}\}$ is  a
$\Sigma^1$-connection from $g$ to $h.$ Then $\{g_1, g_2, g_3, ..., g_{2n+1}, k, l\}$  is  a
$\Sigma^1$-connection from $g$ to $hkl$ in case $g_1 g_2 g_3...g_{2n+1}=h$  and  $\{g_1, g_2, g_3, ..., g_{2n+1}, k^{-1}, l^{-1}\}$   in case $g_1 g_2 g_3...g_{2n+1}=h^{-1}.$  Hence, $hkl\in[g]$ and so we get 
\begin{equation}\label{42}
[\vv_{[g]}, \LL_{1, [g]}, \vv_{[g]}]\subset \vv_{[g]}.
\end{equation}

Finally, consider the first  summand $[\LL_{1, [g]}, \LL_{1, [g]},   \LL_{1, [g]}]$  in (\ref{200}). By Eq. \ref{38}, we have

\begin{eqnarray}\label{43}
\nonumber[\LL_{1, [g]}, \LL_{1, [g]},   \LL_{1, [g]}]&\subset&\sum_{h, k, l\in[g]\cap\Lam^1}[ A_{h^{-1}}\LL_h,  A_{k^{-1}}\LL_k,  A_{l^{-1}}\LL_l]+\\
&~&\sum_{\substack{h, k\in [g]\cap\Lam^1\\ l', l''\in[g]}}\big[ A_{h^{-1}}\LL_h,  A_{k^{-1}}\LL_k, [\LL_{l'}\LL_{l''}, \LL_{{(l'l'')}^{-1}}]\big]\\
\nonumber&+&\sum_{\substack{h, l\in[g]\cap\Lam^1\\k', k''\in[g]}}\big[ A_{h^{-1}}\LL_h, [\LL_{k'}, \LL_{k''}, \LL_{{(k'k'')}^{-1}}], A_{l^{-1}}\LL_l\big]\\
\nonumber&+&\sum_{\substack{h\in[g]\cap\Lam^1\\k', k'', l', l''\in[g]}}\big[ A_{h^{-1}}\LL_h, [\LL_{k'}, \LL_{k''}, \LL_{{(k'k'')}^{-1}}]\\
\nonumber~~&~~&~, [\LL_{l'}, \LL_{l''}, \LL_{{(l'l'')}^{-1}}]\big]\\
\nonumber&+&\sum_{\substack{h', h''\in[g]\\k\in[g]\cap\Lam^1}}\big[[\LL_{h'}, \LL_{h''}, \LL_{{(h'h'')}^{-1}}],  A_{k^{-1}}\LL_k,  A_{l^{-1}}\LL_l\big]\\
\nonumber&+&\sum_{\substack{h', h'', l', l''\in[g]\\k\in[g]\cap\Lam^1}}\big[[\LL_{h'}, \LL_{h''}, \LL_{{(h'h'')}^{-1}}],   A_{k^{-1}}\LL_k\\
\nonumber~~&~~&~,   [\LL_{l'}, \LL_{l''}, \LL_{{(l'l'')}^{-1}}]\big]\\
\nonumber&+&\sum_{\substack{h', h'', k', k''\in[g]\\l\in[g]\cap\Lam^1}}\big[[\LL_{h'}, \LL_{h''}, \LL_{(h'h'')^{-1}}],  [\LL_{k'}, \LL_{k''}, \LL_{{(k'k'')}^{-1}}]\\
\nonumber~~&~~&~,  A_{l^{-1}}\LL_l\big]\\
\nonumber&+&\sum_{h', h'', k', k'', l', l''\in[g]}[[\LL_{h'}, \LL_{h''}, \LL_{(h'h'')^{-1}}],  [\LL_{k'}, \LL_{k''}, \LL_{{(k'k'')}^{-1}}],\\
\nonumber&~&\vspace{10cm}{ [\LL_{l'}, \LL_{l''}, \LL_{{(l'l'')}^{-1}}]]}
\end{eqnarray}

For the first  summand in (\ref{43}), if there exist $h, k, l\in[g]\cap\Lam^1$ such that
$$
[ A_{h^{-1}}\LL_h,  A_{k^{-1}}\LL_k,  A_{l^{-1}}\LL_l]\neq 0,
$$
by  Eqs. (\ref{4}) and (\ref{00}) we have
\begin{eqnarray}
\nonumber[ A_{h^{-1}}\LL_h,  A_{k^{-1}}\LL_k,  A_{l^{-1}}\LL_l]&\subset&[ \LL_{hh^{-1}},  \LL_{kk^{-1}},  A_{l^{-1}}\LL_l]\\
\nonumber&=&A_{l^{-1}}[ \LL_{hh^{-1}},  \LL_{kk^{-1}},  \LL_l]+\rho(\LL_{hh^{-1}},  \LL_{kk^{-1}})A_{l^{-1}}\LL_l\\
\nonumber&\subset&A_{l^{-1}}\LL_l\subset\LL_{1, [g]}.
\end{eqnarray}
For the second  summand in (\ref{43}), if there exist $h, k\in [g]\cap\Lam^1,  l', l''\in[g]$ such that
$$
\big[ A_{h^{-1}}\LL_h,  A_{k^{-1}}\LL_k, [\LL_{l'}\LL_{l''}, \LL_{{(l'l'')}^{-1}}]\big]\neq 0,
$$
by  Eqs. (\ref{4}),  (\ref{00}) and  skew symmetry  we have
\begin{eqnarray}
\nonumber\big[ A_{h^{-1}}\LL_h,  A_{k^{-1}}\LL_k, [\LL_{l'}\LL_{l''}, \LL_{{(l'l'')}^{-1}}]\big]&\subset&\big[ \LL_{hh^{-1}},   [\LL_{l'}\LL_{l''}, \LL_{{(l'l'')}^{-1}}],   A_{k^{-1}}\LL_k\big]\\
\nonumber&=&A_{k^{-1}}\big[\LL_{hh^{-1}},  [\LL_{l'}\LL_{l''}, \LL_{{(l'l'')}^{-1}}],  \LL_k\big]\\
\nonumber&+&\rho\big(\LL_{hh^{-1}},  [\LL_{l'}\LL_{l''}, \LL_{{(l'l'')}^{-1}}]\big)A_{k^{-1}}\LL_k\\
\nonumber&\subset&A_{k^{-1}}\LL_k\subset\LL_{1, [g]}.
\end{eqnarray}
The proof for the rest of summands (exept the last summand)  in (\ref{43})  are similar. For the last summand  in (\ref{43}), if  there exist $h', h'', k', k'', l', l''\in[g]$ such that
$$
\big[[\LL_{h'}, \LL_{h''}, \LL_{(h'h'')^{-1}}],  [\LL_{k'}, \LL_{k''}, \LL_{{(k'k'')}^{-1}}], [\LL_{l'}, \LL_{l''}, \LL_{{(l'l'')}^{-1}}]\big]\neq 0.
$$
By applying identities in Defenition \ref{3L}, we get
\begin{eqnarray}\label{555}
\nonumber\big[[\LL_{h'}, \LL_{h''}, \LL_{(h'h'')^{-1}}],  \LL_1, \LL_1\big]&\subset&\big[[\LL_{h'}, \LL_1, \LL_1], \LL_{h''}, \LL_{(h'h'')^{-1}}\big]\\
&+&\big[[\LL_{h''}, \LL_1, \LL_1], \LL_{(h'h'')^{-1}}, \LL_{h'}\big]\\
\nonumber&+&\big[[\LL_{(h'h'')^{-1}}, \LL_1, \LL_1], \LL_{h'}, \LL_{h''}\big]\\
\nonumber&\subset&[\LL_{h'}, \LL_{h''}, \LL_{(h'h'')^{-1}}]+[\LL_{h''}, \LL_{(h'h'')^{-1}}, \LL_{h'}]\\
\nonumber&+&[\LL_{(h'h'')^{-1}}, \LL_{h'}, \LL_{h''}]\\
\nonumber&\subset&\LL_{1, [g]}.
\end{eqnarray}
Thus, all summands in (\ref{43}) contained in $\LL_{1, [g]}.$ Therefore,
\begin{equation}\label{42.5}
[\LL_{1, [g]}, \LL_{1, [g]}, \LL_{1, [g]}]\subset \LL_{1, [g]}.
\end{equation}
From  Eqs. (\ref{38.5}),  (\ref{41}),  (\ref{42}) and  (\ref{42.5}),  we conclude that  $[I_{[g]}, I_{[g]}, I_{[g]}]\subset I_{[g]}.$

 (2) Observe that
\begin{eqnarray*}
A I_{[g]}&=&\big(A_1\oplus(\bigoplus_{\lam\in\Lam^1}A_\lam)\big)\big(I_{1, [g]}\oplus\vv_{[g]}\big)\\
&=&\big(A_1\oplus(\bigoplus_{\lam\in\Lam^1}A_\lam)\big)\bigg(\big(\sum_{h\in[g]\cap\Lam^1} A_{h^{-1}}\LL_h\big)\\
&+&\big(\sum_{h, k\in[g]}[\LL_h, \LL_k, \LL_{(hk)^{-1}}]\big)\oplus \bigoplus_{h\in[g]}\LL_h\bigg).
\end{eqnarray*}
We discuss it in six cases:

{\bf Case 1.} For the item $A_1(\sum_{h\in[g]\cap\Lam^1} A_{h^{-1}}\LL_h),$ since $\LL$ is an $A$-module, for $h\in[g]\cap\Lam^1$ we have
$$
A_1(A_{h^{-1}}\LL_h)=(A_1 A_{h^{-1}})\LL_h\subset A_{h^{-1}}\LL_h\subset \LL_{1, [g]}. 
$$
Therefore,
$$
A_1(\sum_{h\in[g]\cap\Lam^1} A_{h^{-1}}\LL_h)\subset I_{[g]}.
$$

{\bf Case 2.} Consider the item $A_1\big(\sum_{h, k\in[g]}[\LL_h, \LL_k, \LL_{(hk)^{-1}}]\big).$ By Eq. (\ref{4}), for any $h, k\in[g],$ we have  
 \begin{eqnarray*}
A_1[\LL_h, \LL_k, \LL_{(hk)^{-1}}]&\subset&[\LL_h, \LL_k, A_1\LL_{(hk)^{-1}}]+\rho(\LL_h, \LL_k)(A_1) \LL_{(hk)^{-1}}\\
&\subset&[\LL_h, \LL_k, \LL_{(hk)^{-1}}]+A_{(hk}\LL_{(hk)^{-1}},
\end{eqnarray*}
thanks to $ A_1\LL_{(hk)^{-1}}\subset\LL_{(hk)^{-1}}$ and $\rho(\LL_h, \LL_k)(A_1)\subset A_{hk}.$ Now, if $A_{hk}\neq 0$ (otherwise is trivial), then $hk\in[g]\cap\Lam^1$.  Thus $[\LL_h, \LL_k, \LL_{(hk)^{-1}}]+A_{hk}\LL_{(hk)^{-1}}\subset\LL_{1, [g]}.$ Therefore,
$$
A_1(\sum_{h, k\in[g]}[\LL_h, \LL_k, \LL_{(hk)^{-1}}])\subset I_{[g]}.
$$

{\bf Case 3.} Let us consider the item $A_1(\bigoplus_{h\in[g]}\LL_h).$ Since $\LL$ is an $A$-module, for $h\in[g]$ we have $A_1\LL_h\subset \LL_h\subset\vv_{[g]}.$ Therefore,
$$
A_1(\bigoplus_{h\in[g]}\LL_h)\subset I_{[g]}.
$$

{\bf Case 4.} For the item $(\bigoplus_{\lam\in\Lam^1}A_\lam)(\sum_{h\in[g]\cap\Lam^1} A_{h^{-1}}\LL_h),$ suppose $\lam\in\Lam^1, h\in[g]\cap\Lam^1,$  by associativity  of $A$  we have
$$
A_\lam(A_{h^{-1}}\LL_h)=(A_\lam A_{h^{-1}})\LL_h\subset A_{\lam h^{-1}}\LL_h.
$$
If $\lam h^{-1}\in\Sigma^1$ and $\LL_\lam\neq 0,$  then the family $\{h, \lam  h^{-1}, 1\}$ is a $\Sigma^1$-connection from $h$ to $\Lam.$  That is $\lam\in [g],$ so $ A_{\lam h^{-1}}\LL_h\subset\LL_\lam\subset\vv_{[g]}.$ Therefore
$$ 
(\bigoplus_{\lam\in\Lam^1}A_\lam)(\sum_{h\in[g]\cap\Lam^1} A_{h^{-1}}\LL_h)\subset I_{[g]}. 
$$

{\bf Case 5.} Consider the item $(\bigoplus_{\lam\in\Lam^1}A_\lam)\big(\sum_{h, k\in[g]}[\LL_h, \LL_k, \LL_{(hk)^{-1}}]\big).$  By Eq. (\ref{4}), for $\lam\in\Lam^1, h\in[g]$ we have 
\begin{eqnarray*}
A_\lam[\LL_h, \LL_k, \LL_{(hk)^{-1}}]&\subset&[\LL_h, \LL_k, A_\lam\LL_{(hk)^{-1}}]+\rho(\LL_h, \LL_k)(A_\lam)\LL_{(hk)^{-1}}\\
&\subset&[\LL_h, \LL_k, \LL_{\lam(hk)^{-1}}]+A_{\lam(hk)}\LL_{(hk)^{-1}}.
\end{eqnarray*}
  As in previous case, if $\LL_\lam\neq 0$ and $\LL_{\lam(hk)^{-1}}\neq 0$ we get $\lam, \lam(hk)^{-1}\in\Sigma^1,$  and by Remark \ref{31} we have $\lam\in[g].$ So $[\LL_h, \LL_k, \LL_{\lam(hk)^{-1}}]+A_{\lam(hk)}\LL_{(hk)^{-1}}\subset \vv_{[g]}.$
 Therefore,
$$
(\bigoplus_{\lam\in\Lam^1}A_\lam)\big(\sum_{h, k\in[g]}[\LL_h, \LL_k, \LL_{(hk)^{-1}}]\big)\subset I_{[g]}.
$$

{\bf Case 6.} Finally, consider the item $(\bigoplus_{\lam\in\Lam^1}A_\lam)(\bigoplus_{h\in[g]}\LL_h).$ For $\lam\in\Lam^1$ and $h\in[g]$ we have $A_\lam\LL_h\subset\LL_{\lam h}.$ Using  Remark \ref{31} as in previous case, we can prove $\lam h\in[g].$ Hence, $A_\lam\LL_h\subset\vv_{[g]}.$ Therfore,
$$
(\bigoplus_{\lam\in\Lam^1}A_\lam)(\bigoplus_{h\in[g]}\LL_h)\subset I_{[g]}.
$$ 
Now, summarizing a discussion of above six cases, we get the result.\qed

\begin{pro}\label{12345} Let $[g], [h], [k]\in\Sigma^1/\sim_{\Sigma^1}$ be different from each other,  then 
$$
[I_{[g]}, I_{[h]}, I_{[k]}]=0, ~~~~\hbox{and}~~~~[I_{[g]}, I_{[g]}, I_{[h]}]=0.
$$ 	
\end{pro}
\noindent {\bf Proof.}  We have
\begin{eqnarray}\label{205}
\nonumber[I_{[g]}, I_{[h]}, I_{[k]}]&=&[\LL_{1, [g]}\oplus\vv_{[g]}, \LL_{1, [h]}\oplus\vv_{[h]}, \LL_{1, [k]}\oplus\vv_{[k]}]\\
\nonumber&\subset&[\LL_{1, [g]}, \LL_{1, [h]}, \LL_{1, [k]}]+[\LL_{1, [g]}, \LL_{1, [h]}, \vv_{[k]}]+[\LL_{1, [g]}, \vv_{[h]}, \LL_{1, [k]}]\\
\nonumber&+&[\LL_{1, [g]}, \vv_{[h]}, \vv_{[k]}]+[\vv_{[g]}, \LL_{1, [h]}, \LL_{1, [k]}]+[\vv_{[g]}, \LL_{1, [h]}, \vv_{[k]}]\\
&+&[\vv_{[g]}, \vv_{[h]}, \LL_{1, [k]}]+[\vv_{[g]}, \vv_{[h]}, \vv_{[k]}]. 
\end{eqnarray}
Let us consider the last  item $[\vv_{[g]}, \vv_{[h]}, \vv_{[k]}]$ in Eq. (\ref{205}).  Suppose that there exist $g_1\in [g], h_1\in [h]$ and $k_1\in [k]$  such that
$[\LL_{g_1}, \LL_{h_1}, \LL_{k_1}]\neq 0.$ By definition of grading, $g_1h_1k_1\in\Sigma^1.$ Since  $g_1\in[g]$ and $g_1h_1k_1\in\Sigma^1,$ we get $g\sim_{\Sigma^1} g_1h_1k_1.$ Similarly, one can get  $h\sim_{\Sigma^1} g_1h_1k_1.$ Now, Proposition \ref{rela}  implies that $[g]=[h]$ a contradiction. Therefore,
\begin{equation}\label{430}
[\vv_{[g]}, \vv_{[h]}, \vv_{[k]}]=0.
\end{equation}

Now, we consider the item $[\LL_{1, [g]}, \vv_{[h]}, \vv_{[k]}]$ in Eq. (\ref{205}). We have
\begin{eqnarray}\label{431}
\nonumber[\LL_{1, [g]}, \vv_{[h]}, \vv_{[k]}] &=&\bigg[\big(\sum_{g'\in[g]\cap\Lam^1} A_{g'^{-1}}\LL_{g'}\big)+\big(\sum_{h', k'\in[g]}[\LL_{h'}, \LL_{k'}, \LL_{(h'k')^{-1}}]\big), \\
&~&\bigoplus_{h_1\in[h]} \LL_{h_1}, \bigoplus_{k_1\in[k]}\LL_{k_1}\bigg]\\
\nonumber&\subset&\big[\sum_{g'\in[g]\cap\Lam^1} A_{g'^{-1}}\LL_{g'}, \bigoplus_{h_1\in[h]} \LL_{h_1}, \bigoplus_{k_1\in[k]}\LL_{k_1}\big]\\
 \nonumber&+&\big[\sum_{h', k'\in[g]}[\LL_{h'}, \LL_{k'}, \LL_{(h'k')^{-1}}], \bigoplus_{h_1\in[h]} \LL_{h_1}, \bigoplus_{k_1\in[k]}\LL_{k_1}\big].
\end{eqnarray}
For the first summand in (\ref{431}),  suppose  there exist $g'\in[g]\cap\Lam^1, h_1\in[h]$ and $k_1\in[k]$ such that $[ A_{g'^{-1}}\LL_{g'}, \LL_{h_1}, \LL_{k_1}]\neq 0.$ By Eq. (\ref{4}),
\begin{equation*}
[ A_{g'^{-1}}\LL_{g'}, \LL_{h_1}, \LL_{k_1}]=A_{g'^{-1}}[\LL_{h_1}, \LL_{k_1}, \LL_{g'}]+\rho( \LL_{h_1}, \LL_{k_1})A_{g'^{-1}}\LL_{g'}.
 \end{equation*}
Taking into account  Eq.(\ref{430}), we get $ [\LL_{h_1}, \LL_{k_1}, \LL_{g'}]=0.$ If $\rho( \LL_{h_1}, \LL_{k_1})A_{g'^{-1}}\LL_{g'}\neq 0,$ then $A_{h_1k_1{g'}^{-1}}\neq 0$ and $h_1k_1{g'}^{-1}\in\Lam^1.$
 We take the family $\{g', k_1^{-1}, h_1k_1{g'}^{-1}\}$ as a $\Sigma^1$-connection from $g'$ to $h_1,$ and so $[g]=[h]$ wich is a contradiction. That is   $\rho( \LL_{h_1}, \LL_{k_1})A_{g'^{-1}}\LL_{g'}=0.$ Hence, 
\begin{equation}\label{431.01}
[ A_{g'^{-1}}\LL_{g'}, \LL_{h_1}, \LL_{k_1}]=0,
 \end{equation}
and so
\begin{equation}\label{431.1}
[\sum_{g'\in[g]\cap\Lam^1} A_{g'^{-1}}\LL_{g'}, \bigoplus_{h_1\in[h]} \LL_{h_1}, \bigoplus_{k_1\in[k]}\LL_{k_1}]=0.
 \end{equation}
Next, consider the second  summand in (\ref{431}). For $h', k'\in[g], h_1\in[h]$ and $k_1\in[k],$
 by fundamental identity and Eq. (\ref{430}), we get
\begin{eqnarray*}
\big[[\LL_{h'}, \LL_{k'}, \LL_{(h'k')^{-1}}],  \LL_{h_1}, \LL_{k_1}\big]&\subset& \big[[\LL_{h'},  \LL_{h_1}, \LL_{k_1}], \LL_{k'}, \LL_{(h'k')^{-1}}\big]\\
&+& \big[[\LL_{k'},  \LL_{h_1}, \LL_{k_1}], \LL_{(h'k')^{-1}}, \LL_{h'}\big]\\
&+& \big[[\LL_{(h'k')^{-1}},  \LL_{h_1}, \LL_{k_1}], \LL_{h'}, \LL_{k'}\big]\\
&=&0,
\end{eqnarray*}
and so
\begin{equation}\label{431.2}
\big[\sum_{h', k'\in[g]}[\LL_{h'}, \LL_{k'}, \LL_{(h'k')^{-1}}], \bigoplus_{h_1\in[h]} \LL_{h_1}, \bigoplus_{k_1\in[k]}\LL_{k_1}\big]=0.
\end{equation}
 From Eqs. (\ref{431.1}) and (\ref{431.2}), we have  $[\LL_{1, [g]}, \vv_{[h]}, \vv_{[k]}]=0.$
 By  the skew symmetry of  trilinear map, we also get
\begin{equation}\label{432}
[\LL_{1, [g]}, \vv_{[h]}, \vv_{[k]}]=[\vv_{[g]}, \LL_{1, [h]}, \vv_{[k]}]=[\vv_{[g]}, \vv_{[h]}, \LL_{1, [k]}]=0.
\end{equation}

Next, we consider the summand $[\LL_{1, [g]}, \LL_{1, [h]}, \vv_{[k]}]$ in Eq. (\ref{205}). We have
 \begin{eqnarray}\label{404}
 \nonumber[\LL_{1, [g]}, \LL_{1, [h]}, \vv_{[k]}] &=&\bigg[\big(\sum_{g'\in[g]\cap\Lam^1} A_{{g'}^{-1}}\LL_{g'}\big)+\big(\sum_{{g_1, g_2}\in[g]}[\LL_{g_1}, \LL_{g_2}, \LL_{(g_1g_2)^{-1}}]\big),\\
\nonumber&~& \big(\sum_{h'\in[h]\cap\Lam^1} A_{{h'}^{-1}}\LL_{h'}\big)+\big(\sum_{h_1, h_2\in[h]}[\LL_{h_1}, \LL_{h_2}, \LL_{(h_1h_2)^{-1}}]\big), \\
&~& \bigoplus_{k_1\in[k]}\LL_{k_1}\bigg].
 \end{eqnarray}
 The above statement includs four items which we consider in the following. First, consider the item $[\sum_{g'\in[g]\cap\Lam^1} A_{{g'}^{-1}}\LL_{g'}, \sum_{h'\in[h]\cap\Lam^1} A_{{h'}^{-1}}\LL_{h'}, \bigoplus_{k_1\in[k]}\LL_{k_1}]$ in (\ref{404}). For  $g'\in[g]\cap\Lam^1, h'\in[h]\cap\Lam^1$ and $k_1\in[k],$ by  using Eqs. (\ref{430}) and (\ref{431.01}) we have
 
\begin{eqnarray*}
 [ A_{{g'}^{-1}}\LL_{g'},   A_{{h'}^{-1}}\LL_{h'}, \LL_{k_1}]&=& [ \LL_{k_1},   A_{{h'}^{-1}}\LL_{h'},  A_{{g'}^{-1}}\LL_{g'}]\\
&=&A_{{g'}^{-1}}[\LL_{k_1},   A_{{h'}^{-1}}\LL_{h'}, \LL_{g'}]+\rho(\LL_{k_1},   A_{{h'}^{-1}}\LL_{h'})A_{{g'}^{-1}}\LL_{g'}\\ 
&=& A_{{g'}^{-1}}\big( A_{{h'}^{-1}}[\LL_{k_1},   \LL_{h'}, \LL_{g'}]+\rho(\LL_{k_1},  \LL_{g'}) A_{{h'}^{-1}}\LL_{h'}\big)\\
&+&A_{{h'}^{-1}}\rho(\LL_{k_1},   \LL_{h'})A_{{g'}^{-1}}\LL_{g'}\\
&=& A_{{g'}^{-1}} A_{{h'}^{-1}}[\LL_{k_1},   \LL_{h'}, \LL_{g'}]+A_{{g'}^{-1}}\rho(\LL_{k_1},  \LL_{g'}) A_{{h'}^{-1}}\LL_{h'}\\
&+&A_{{h'}^{-1}}\rho(\LL_{k_1},   \LL_{h'})A_{{g'}^{-1}}\LL_{g'}\\
&=&0.
\end{eqnarray*}
 Therefore,  
\begin{equation}\label{404.5}
[\sum_{g'\in[g]\cap\Lam^1} A_{{g'}^{-1}}\LL_{g'}, \sum_{h'\in[h]\cap\Lam^1} A_{{h'}^{-1}}\LL_{h'}, \bigoplus_{k_1\in[k]}\LL_{k_1}]=0.
\end{equation}
Second, consider the item $$\big[\sum_{g'\in[g]\cap\Lam^1} A_{{g'}^{-1}}\LL_{g'}, \sum_{h_1, h_2\in[h]}[\LL_{h_1}, \LL_{h_2}, \LL_{(h_1h_2)^{-1}}], \bigoplus_{k_1\in[k]}\LL_{k_1}\big],$$  in Eq. (\ref{404}).  For  $g'\in[g]\cap\Lam^1, h_1, h_2\in[h]$  and $k_1\in[k],$  again  by  using Eqs. (\ref{430}) and (\ref{431.01}) we have
\begin{eqnarray*}
\big[ A_{{g'}^{-1}}\LL_{g'},   [\LL_{h_1}, \LL_{h_2}, \LL_{(h_1h_2)^{-1}}], \LL_{k_1}\big]&=&  \big[[\LL_{h_1}, \LL_{h_2}, \LL_{(h_1h_2)^{-1}}], \LL_{k_1}, A_{{g'}^{-1}}\LL_{g'}\big]\\
&=&A_{{g'}^{-1}} \big[[\LL_{h_1}, \LL_{h_2}, \LL_{(h_1h_2)^{-1}}], \LL_{k_1}, \LL_{g'}\big]\\
&+&\rho\big([\LL_{h_1}, \LL_{h_2}, \LL_{(h_1h_2)^{-1}}], \LL_{k_1}\big)A_{{g'}^{-1}}\LL_{g'}\\
&=&0.
\end{eqnarray*}
Hence,
\begin{equation}\label{406}
\big[\sum_{g'\in[g]\cap\Lam^1} A_{{g'}^{-1}}\LL_{g'}, \sum_{h_1, h_2\in[h]}[\LL_{h_1}, \LL_{h_2}, \LL_{(h_1h_2)^{-1}}], \bigoplus_{k_1\in[k]}\LL_{k_1}\big]=0.
\end{equation}
 By skew symmetry, we also have  
\begin{equation}\label{407}
\big[\sum_{{g_1, g_2}\in[g]}[\LL_{g_1}, \LL_{g_2}, \LL_{(g_1g_2)^{-1}}], \sum_{h'\in[h]\cap\Lam^1} A_{{h'}^{-1}}\LL_{h'},  \bigoplus_{k_1\in[k]}\LL_{k_1}\big]=0.
\end{equation}
Now, consider the forth item in Eq. (\ref{404}). For $g_1, g_2\in[g], h_1, h_2\in[h]$ and $k_1\in[k],$ using Eq. (\ref{430}) we have 
$$
big[[\LL_{g_1}, \LL_{g_2}, \LL_{(g_1g_2)^{-1}}], [\LL_{h_1}, \LL_{h_2}, \LL_{(h_1h_2)^{-1}}], \LL_{k_1}]\subset[\vv_{[g]}, \vv_{[h]}, \vv_{[k]}\big]=0.
$$
 So we get
\begin{equation}\label{409}
\big[\sum_{{g_1, g_2}\in[g]}[\LL_{g_1}, \LL_{g_2}, \LL_{(g_1g_2)^{-1}}], \sum_{h_1, h_2\in[h]}[\LL_{h_1}, \LL_{h_2}, \LL_{(h_1h_2)^{-1}}],  \bigoplus_{k_1\in[k]}\LL_{k_1}\big]=0.
\end{equation}
Now, from Eqs. (\ref{404.5})-(\ref{409}), we get
\begin{equation}\label{411}
[\LL_{1, [g]}, \LL_{1, [h]}, \vv_{[k]}]=0.
\end{equation}
 By skew symmetry, we also get
\begin{equation}\label{412.5}
[\LL_{1, [g]}, \vv_{[h]}, \LL_{1, [k]}]=[\vv_{[g]}, \LL_{1, [h]}, \LL_{1, [k]}]=[\LL_{1, [g]}, \LL_{1, [h]}, \vv_{[k]}]=0.
\end{equation}

Finally, consider the first item  in Eq. (\ref{205}). We have
\begin{eqnarray}\label{412}
 \nonumber[\LL_{1, [g]}, \LL_{1, [h]}, \LL_{1, [k]}] &=&\bigg[\big(\sum_{g'\in[g]\cap\Lam^1} A_{{g'}^{-1}}\LL_{g'}\big)+\big(\sum_{{g_1, g_2}\in[g]}[\LL_{g_1}, \LL_{g_2}, \LL_{(g_1g_2)^{-1}}]\big),\\
\nonumber~~&~~&~ \LL_{1, [h]}, \LL_{1, [k]}\bigg]\\
\nonumber&\subset& \big[\sum_{g'\in[g]\cap\Lam^1} A_{{g'}^{-1}}\LL_{g'}, \LL_{1, [h]}, \LL_{1, [k]}\big]\\
&+&\big[\sum_{{g_1, g_2}\in[g]}[\LL_{g_1}, \LL_{g_2}, \LL_{(g_1g_2)^{-1}}], \LL_{1, [h]}, \LL_{1, [k]}\big].
 \end{eqnarray}
By a similar argument as in Eq. (\ref{407}), the first summand in (\ref{412}) is zero. For the second summand (\ref{412}), suppose $g_1, g_2\in[g]$ we have
\begin{eqnarray*}
 \big[[\LL_{g_1}, \LL_{g_2}, \LL_{(g_1g_2)^{-1}}], \LL_{1, [h]}, \LL_{1, [k]}\big]&\subset&\big[[\LL_{g_1}, \LL_{1, [h]}, \LL_{1, [k]}], \LL_{g_2}, \LL_{(g_1g_2)^{-1}}\big]\\
&+&\big[[\LL_{g_2}, \LL_{1, [h]}, \LL_{1, [k]}], \LL_{(g_1g_2)^{-1}}, \LL_{g_1}\big]\\
&+&\big[[ \LL_{(g_1g_2)^{-1}}, \LL_{1, [h]}, \LL_{1, [k]}], \LL_{g_1}, \LL_{g_2}\big].
 \end{eqnarray*}
All of the above three snmmands are zero, thanks to Eq.  (\ref{412.5}). Therefore,
\begin{equation}\label{413}
[\LL_{1, [g]}, \LL_{1, [h]}, \LL_{1, [k]}]=0.
\end{equation}

From Eqs. (\ref{430}) ,  (\ref{432}),   (\ref{412.5}) and (\ref{413}), we get
$$
[I_{[g]}, I_{[h]}, I_{[k]}]=0.
$$
By a similar argument as above,  one can prove $[I_{[g]}, I_{[g]}, I_{[h]}]=0.$\qed

\begin{thm}\label{ideal} The following assertions hold
\begin{itemize}
\item[(1)] For any $[g]\in\Sigma^1/\sim_{\Sigma^1},$ the linear space
$$
I_{[g]} =\LL_{1, [g]}\oplus\vv_{[g]},
$$
associated to $[g]$ is a graded  ideal of $(\LL, A).$
\item[(2)]  If $(\LL, A)$ is gr-simple, then there exists a $\Sigma^1$-connection from
$g$ to $h$ for any $g, h\in\Sigma^1,$ and
$$
\LL_1=\sum_{g\in\Sigma^1\cap\Lam^1} A_{g^{-1}}\LL_g
+\sum_{h, k\in\Sigma^1}[\LL_h, \LL_k, \LL_{(hk)^{-1}}].
$$
\end{itemize}
\end{thm}
\noindent {\bf Proof.}
(1) We are going to check $[I_{[g]},\LL, \LL]\subset I_{[g]}.$ We have
\begin{eqnarray}\label{421}
 \nonumber[I_{[g]},\LL, \LL] &=&\big[\LL_{1, [g]}\oplus\vv_{[g]}, \LL_1\oplus(\bigoplus_{h\in \Sigma^1}\LL_h), \LL_1\oplus(\bigoplus_{k\in \Sigma^1}\LL_k)\big]\\
\nonumber&\subset& [\LL_{1, [g]}, \LL_1, \LL_1]+[\LL_{1, [g]}, \LL_1, \bigoplus_{k\in \Sigma^1}\LL_k]+[\LL_{1, [g]}, \bigoplus_{h\in \Sigma^1}\LL_h, \LL_1]\\
\nonumber&+&[\LL_{1, [g]}, \bigoplus_{h\in \Sigma^1}\LL_h, \bigoplus_{k\in \Sigma^1}\LL_k]
+[\vv_{[g]}, \LL_1, \LL_1]+[\vv_{[g]}, \LL_1, \bigoplus_{k\in \Sigma^1}\LL_k]\\
&+&[\vv_{[g]}, \bigoplus_{h\in \Sigma^1}\LL_h, \LL_1]+[\vv_{[g]}, \bigoplus_{h\in \Sigma^1}\LL_h, \bigoplus_{k\in \Sigma^1}\LL_k].
 \end{eqnarray}
Let us consider the first summan in Eq. (\ref{421}), we have
\begin{eqnarray}\label{422}
 \nonumber [\LL_{1, [g]}, \LL_1, \LL_1]&=&\big[\sum_{h\in[g]\cap\Lam^1} A_{h^{-1}}\LL_h+\sum_{h, k\in[g]}[\LL_h, \LL_k, \LL_{(hk)^{-1}}], \LL_1, \LL_1\big]\\
&\subset&[\sum_{h\in[g]\cap\Lam^1} A_{h^{-1}}\LL_h, \LL_1, \LL_1]+\big[\sum_{h, k\in[g]}[\LL_h, \LL_k, \LL_{(hk)^{-1}}], \LL_1, \LL_1\big].
 \end{eqnarray}
Suppose $h\in[g]\cap\Lam^1,$  by Eq. (\ref{4}),
$$
[ A_{h^{-1}}\LL_h, \LL_1, \LL_1]=A_{h^{-1}}[\LL_1, \LL_1, \LL_h]+\rho(\LL_1, \LL_1)A_{h^{-1}}\LL_h\subset\LL_{1, [g]}.
$$
Now, if  $h, k\in[g]$ then $\big[[\LL_h, \LL_k, \LL_{(hk)^{-1}}], \LL_1, \LL_1\big]\subset\LL_{1, [g]},$ thanks to Eq. (\ref{555}). Taking into account Eq. (\ref{422}), we get 
$$
[\LL_{1, [g]}, \LL_1, \LL_1]\subset\LL_{1, [g]}.
$$
Next, by Proposition \ref{subalg}(1), for the rest of all summands in Eq. (\ref{422}), we get $[I_{[g]},\LL, \LL]\subset I_{[g]}.$ So $I_{[g]}$ is a $3$-Lie ideal of $\LL.$ By  Proposition \ref{subalg}(2), we also have $A I_{[g]}\subset I_{[g]},$ that is $I_{[g]}$ is an $A$-module. Finally,  by Eq. (\ref{4}) we have
$$
\rho(I_{[g]}, I_{[g]})(A)\LL\subset[I_{[g]}, I_{[g]}, A\LL ]+A[I_{[g]}, I_{[g]}, \LL]\subset[I_{[g]}, I_{[g]}, \LL ]+I_{[g]}\subset I_{[g]}.
$$
By construction of $I_{[g]},$ it is  a graded ideal of $(\LL, A).$\\

(2)  The gr-simplicity of $(\LL, A)$ implies that $0\neq I_{[g]}\in\{\LL, \ker\rho\}$ for any $g\in\Sigma^1.$  If $I_{[g]}=\LL$ for some $g\in\Sigma^1,$  then $[g]=\Sigma^1.$ Otherwise, if $I_{[g]}=\ker\rho$ for all $g\in\Sigma^1$ we have $[g]=[h]$ for any $h\in\Sigma^1$  and again $\Sigma^1=[g].$ We  conclude that  all the  elements of the $G$-support $\Sigma^1$ are  $\Sigma^1$-connected . Moreover, clearly
$$
\LL_1=\sum_{g\in\Sigma^1\cap\Lam^1} A_{g^{-1}}\LL_g
+\sum_{h, k\in\Sigma^1}[\LL_h, \LL_k, \LL_{(hk)^{-1}}].
$$
\qed

\begin{thm}\label{main1} Let  $(\LL, A)$ be a graded  $3$-Lie-Rinehart  algebra.  For a vector space complement  $\uu$ of   
$
\sum_{g\in\Sigma^1\cap\Lam^1} A_{g^{-1}}\LL_g
+\sum_{h, k\in\Sigma^1}[\LL_h, \LL_k, \LL_{(hk)^{-1}}],
$
in $\LL_1,$  we have
$$
\LL=\uu\oplus\sum_{[g]\in\Sigma^1/\sim_{\Sigma^1}}I_{[g]},
$$
where any $I_{[g]}$ is one of the graded ideals of $(\L, A)$ described in Theorem \ref{ideal}-(1). Furthermore, $ [I_{[g]}, I_{[h]}, I_{[k]}]=0$ where $[g], [h], [k]\in\Sigma^1/\sim_{\Sigma^1}$ be different from each other.
\end{thm}
\noindent {\bf Proof.} Each $I_{[g]}$ is well defined and  by
Theorem \ref{ideal}-(1), a graded ideal of $(\LL, A).$ It is
clear that
$$
\LL=\LL_1\oplus(\bigoplus_{g\in\Sigma^1}\LL_g)=\uu\oplus\sum_{[g]\in\Sigma^1/\sim_{\Sigma^1}}I_{[g]},
$$
where $\u$ is a linear  space complement   of   
$$
\sum_{g\in\Sigma^1\cap\Lam^1} A_{g^{-1}}\LL_g
+\sum_{h, k\in\Sigma^1}[\LL_h, \LL_k, \LL_{(hk)^{-1}}],
$$
in $\LL_1.$ 
By  Proposition \ref{12345}, we also have   $ [I_{[g]}, I_{[h]}, I_{[k]}]=0,$ where $[g], [h], [k]\in\Sigma^1/\sim_{\Sigma^1}$ be different from each other. \qed

\begin{cor}\label{2.17} If $Z_{\rho}(\LL)=\{0\}$ and 
$$
\LL_1=\sum_{g\in\Sigma^1\cap\Lam^1} A_{g^{-1}}\LL_g
+\sum_{h, k\in\Sigma^1}[\LL_h, \LL_k, \LL_{(hk)^{-1}}],
$$	
then $\LL$ is the direct sum of the graded ideals
given in Theorem \ref{ideal}-(1),
$$
\LL=\bigoplus_{[g]\in\Sigma^1/\sim_{\Sigma^1}}I_{[g]},
$$	
Moreover, $ [I_{[g]}, I_{[h]}, I_{[k]}]=0$ where $[g], [h], [k]\in\Sigma^1/\sim_{\Sigma^1}$ be different from each other.
\end{cor}\label{}
\noindent {\bf Proof.}  Since $
\LL_1=\sum_{g\in\Sigma^1\cap\Lam^1} A_{g^{-1}}\LL_g
+\sum_{h, k\in\Sigma^1}[\LL_h, \LL_k, \LL_{(hk)^{-1}}],
$
we get
$$
\LL=\sum_{[g]\in\Sigma^1/\sim_{\Sigma^1}}I_{[g]},
$$
where $I_{[g]}$ is one of the graded ideals of $(\L, A)$ described in Theorem \ref{ideal}-(1) satisfying in Proposition \ref{12345}.
For the direct character, suppose there exists  $x\in I_{[g]}\cap\sum_{[h]\in\Sigma^1/\sim_{\Sigma^1}}I_{[h]}$ such that   $[g]\neq [h].$ The fact $[I_{[g]}, I_{[g]}, I_{[h]}]=0$ with $[g]\neq [h]$ and $x\in I_{[g]},$ implies that
$$
[x,\sum_{[g]\in\Sigma^1/\sim_{\Sigma^1}}I_{[g]}, \LL ]=0.
$$
We also have $[x, I_{[g]}, \LL ]=0,$ thanks to $x\in\sum_{[h]\in\Sigma^1/\sim_{\Sigma^1}}I_{[h]}$ with $[g]\neq [h]$ and the same above fact.Therefore, $[x, \LL, \LL]=0.$ Next, by Eq. (\ref{4}), we have $\rho(x, \LL)=0.$ Thus, we get $x\in Z_{\rho}(\LL)=\{0\}.$\qed

\section{ Connections in $\Lam^1$ and decomposition of $A.$}
\setcounter{equation}{0}\

Let $(\LL, A)$  be a $G$-graded $3$-Lie-Rinehartd algebra (see Definition \ref{maindef}). In this section we begin by introducing the so called  connection among of the
elements  in the $G$-support $\Lam^1$ for an assopciative and commutative algebra $A$ associated to $(\LL, A).$ Recall that $A$ admits a $G$-grading as
$$
A=A_1\oplus(\bigoplus_{\lam\in \Lambda^1}A_\lam),
$$
where $\Lam^1=\{\lam\in G\setminus\{1\}~:~A_\lam\neq 0\},$ is the $G$-support of grading. We will consider the sets $\Sigma^{\pm}$ and $\Lam^{\pm}$ as in Section 3. Note that, the proof of some results in this section are similar to the one for graded Lie-Rinehard algebra (see \cite{BCNS}), we will omit them.

\begin{DEF}\label{conn2}
Let $\lam, ~\mu\in \Lam^1,$ we say that {\em $\lam$
is $\Lam^1$-connected to $\mu$} and denoted by $\lam\approx_{\Lam^1}\mu,$  if  there exists a family
$\{\lam_1, \lam_2, \lam_3, ..., \lam_n\}\subset\Sigma\cup\Lam\cup\{1\},$
such that	satisfying the following conditions;\\
\begin{itemize}	
\item[(1)] $\lam_1=\lam.$
	
\item[(2)] $\lam_1\lam_2\in\Lam,\\
		\lam_1 \lam_2 \lam_3\in\Lam,\\
		...\\
		\lam_1 \lam_2 \lam_3...+\lam_{n-1}\in\Lam.$
		
\item[(3)] $\lam_1 \lam_2 \lam_3...\lam_n\in\{\mu, \mu^{-1}\}.$
	\end{itemize}
	The family $\{\lam_1, \lam_2, \lam_3, ..., \lam_n\}$ is called a
	{\em $\Lam^1$-connection} from $\lam$ to $\mu.$
\end{DEF}
 
The next result shows that the $\Lam^1$-connection relation is an equivalence relation (see Proposition 3.2 in \cite{BCNS}).
\begin{pro}\label{rela2}
The relation $\approx_{\Lam^1}$ in $\Lam^1$ defined by
$$
\lam\approx_{\Lam^1}\mu~\hbox{ if and only if}~ \lam~\hbox{ is $\Lam^1-$connected to}~ \mu,
$$
is an equivalence relation.
\end{pro}

\begin{rem}\label{lem}
Let $\lam, \mu\in\Lam^1$ such that $\lam\approx_{\Lam^1}\mu.$ If $\lam\eta\in\Lam^1,$ for $\eta\in\Sigma\cup\Lam$ then $\lam\approx_{\Lam^1}\mu\eta.$ Considering the connection $\{\mu, \eta\}$ we get $\mu\approx_{\Lam^1}\mu\eta$ and by transitivity $\lam\approx_{\Lam^1}\mu\eta.$
\end{rem}

By the  Proposition \ref{rela2}, we can consider the equivalence relation
in $\Lam^1$ by the $\Lam^1$-connection relation $\approx_{\Lam^1}$ in $\Lam^1.$ So we denote
by
$$
\Lam^1/\approx_{\Lam^1}  :=\{[\lam] : \lam\in\Lam^1\},
$$
where $[\lam]$ denotes  the set of  elements of $\Lam^1,$  which are
$\Lam^1$-connected to $\lam.$ 

Our next goal in this section is to associate an adequate graded  ideal
$\aa_{[\lam]}$ of $A$ to any $[\lam]\in\Lam^1/\approx_{\Lam^1}.$ For a fixed  $\lam\in\Lam^1,$ we define
\begin{eqnarray}\label{119}
A_{1, [\lam]} :=\big(\sum_{\mu\in[\lam]} A_{\mu^{-1}}A_\mu\big)+\big(\sum_{h, k\in {[\lam]\cap\Sigma^1}}\rho(\LL_h, \LL_k) A_{(hk)^{-1}}\big)\subset A_1.
\end{eqnarray}
Next, we define
\begin{equation}\label{48'}
A_{[\lam]} :=\bigoplus_{\mu\in[\lam]}A_\mu.
\end{equation}
Finally, we denote by $\aa_{[\lam]}$ the direct sum of the two graded
subspaces above, that is,
\begin{equation}\label{48''}
	\aa_{[\lam]} :=A_{1, [\lam]}\oplus A_{[\lam]}.
\end{equation}
The detail proofs of the following properties of $A$ can be found in  \cite{BCNS};

\begin{pro}\label{4.7} For any  $[\lam]\in\Lam^1/\approx_{\Lam^1},$ we have $\aa_{[\lam]} \aa_{[\lam]}\subset\aa_{[\lam]}.$
\end{pro}

\begin{pro}\label{4.8} For any $\lam, \mu\in\Lam^1,$  if $[\lam]\neq[\mu]$ then $\aa_{[\lam]}\aa_{[\mu]}=0.$
\end{pro}

We recall that a $G$-graded subspace $I$ of a commutative and associative algebra $A$ is called an ideal of $A$ if $A I \subset I.$ We say that $A$ is gr-simple if $A A\neq 0$ and it contains no proper ideals.

\begin{thm}\label{main1'} Let $A$ be a commutative and associative  algebra associated to a  graded $3$-Lie-Rinehart  algebra $(\LL, A).$ Then the following assertions hold.
\begin{itemize}
\item[(1)] For any $[\lam]\in\Lam^1/\approx_{\Lam^1},$ the linear  subspace
$$
\aa_{[\lam]} =A_{1, [\lam]}\oplus A_{[\lam]},
$$
of  algebra $A$ associated to $[\lam]$ is a graded ideal of $A.$
		
\item[(2)]  If $A$  is gr-simple then all elements of $\Lam^1$ are $\Lam^1$-connected. Furthermore,
$$
A_1=\sum_{\mu\in\Lam^1} A_{\mu^{-1}}A_\mu+\sum_{h, k\in {\Lam^1\cap\Sigma^1}}\rho(\LL_h, \LL_k) A_{(hk)^{-1}}.
$$
\end{itemize}
\end{thm}
\noindent {\bf Proof.} The proof is similar to the one in \cite{BCNS}, Theorem  3.6 for a graded Lie-Rinehart algebra.\qed

\begin{thm}\label{main2'} Let $A$ be a commutative and associative  algebra associated to a  graded $3$-Lie-Rinehart  algebra $(\LL, A).$ Then 
$$
A=\v+\sum_{[\lam]\in\Lam^1/\approx_{\Lam^1}}\aa_{[\lam]},
$$
where $\v$ is a linear complement of 
$$
\sum_{\mu\in\Lam^1} A_{\mu^{-1}}A_\mu+\sum_{h, k\in {\Lam^1\cap\Sigma^1}}\rho(\LL_h, \LL_k) A_{(hk)^{-1}},
$$ 
 in $A_1$ and any $\aa_{[\lam]}$ is one of the graded  ideals of $A$ described in Theorem \ref{main1'}-(1). Forthermore  $\aa_{[\lam]} \aa_{[\mu]}=0,$ when $[\lam]\neq [\mu].$
\end{thm}
\noindent {\bf Proof.} The proof is similar to the one in \cite{BCNS}, Theorem  3.7 for a  graded Lie-Rinehart algebra.\qed

Recall that,  denote by 
$$
Ann(A) := \{a \in A : a A = 0\},~~~~~~\hbox{and~}~~~~~Ann_{\LL}(A):= \{x\in\LL :  Ax = 0\},
$$
 the annihilator of the commutative and
associative algebra $A$ and the annihilator of $A$ in $\LL.$

\begin{cor}\label{4.19} Let $(\LL, A)$ be a graded $3$-Lie-Rinehart  algebra. If $Ann(A) = 0$ and
$$
A_1=\sum_{\mu\in\Lam^1} A_{\mu^{-1}}A_\mu+\sum_{h, k\in {\Lam^1\cap\Sigma^1}}\rho(\LL_h, \LL_k) A_{(hk)^{-1}},
$$
then $A$ is the direct sum of the graded  ideals given in Theorem \ref{main1'}-(1),
$$
A=\bigoplus_{[\lam]\in\Lam^1/\approx_{\Lam^1}}\aa_{[\lam]}.
$$
Furthermore, $\aa_{[\lam]} \aa_{[\mu]}=0,$ when $[\lam]\neq [\mu].$
\end{cor}
\noindent {\bf Proof.} This can be proved analogously to Corollary 3.8 in \cite{BCNS}.\qed\\

In the following, we will discuss the relation between the decompositions of $L$ and $A$
of a graded $3$-Lie Rinehart algebra $(\LL, A).$

\begin{DEF}  A graded $3$-Lie-Rinehart  algebra $(\LL, A)$ is tight if $Z_{\rho}(\LL)= 0,~~~~~~Ann(A)=0=Ann_{\LL}(A),~~AA = A,~~ A\LL = \LL$ and
\begin{eqnarray*}
\LL_1=\sum_{g\in\Sigma^1\cap\Lam^1} A_{g^{-1}}\LL_g
+\sum_{h, k\in\Sigma^1}[\LL_h, \LL_k, \LL_{(hk)^{-1}}],\\
A_1=\sum_{\mu\in\Lam^1} A_{\mu^{-1}}A_\mu+\sum_{h, k\in {\Lam^1\cap\Sigma^1}}\rho(\LL_h, \LL_k) A_{(hk)^{-1}}.
\end{eqnarray*}
\end{DEF} 

\begin{rem}\label{rem} If $(\LL, A)$ is a tight  graded $3$-Lie-Rinehart  algebra then it follows from   Corollaries \ref{2.17} and  \ref{4.19} that 
$$
\LL=\bigoplus_{[g]\in\Sigma^1/\sim_{\Sigma^1}}I_{[g]},~~ A=\bigoplus_{[\lam]\in\Lam^1/\approx_{\Lam^1}}\aa_{[\lam]},
$$	
with any $I_{[g]}$ a graed ideal of $\LL$ satisfying $[I_{[g]}, I_{[h]}]=0$ if $[g]\neq[h]$  and any $\aa_{[\lam]}$ a graded ideal of $A$  satisfying  $\aa_{[\lam]} \aa_{[\mu]}=0,$ when $[\lam]\neq [\mu].$
\end{rem}

\begin{pro}\label{5.10}   Let $(\LL, A)$ be a tight graded $3$-Lie-Rinehart  algebra. Then for any $[g]\in\Sigma^1/\sim_{\Sigma^1}$
there exists a unique $[\lam]\in\Lam^1/\approx_{\Lam^1}$ such that $ \aa_{[\lam]} I_{[g]} \neq 0.$ 
\end{pro}
\noindent {\bf Proof.} At first, we are going  to prove the existence. We claim $AI_{[g]}\neq 0$  for any $[g]\in\Sigma^1/\sim_{\Sigma^1}.$ Indeed, if $AI_{[g]}= 0$ for some $[g]\in\Sigma^1/\sim_{\Sigma^1},$ then by the fact $I_{[g]}$ is a graded  ideal of $\LL$ we have 
$$
[I_{[g]}, A\LL, A\LL]=\big[I_{[g]}, \bigoplus_{[h]\in\Sigma^1/\sim_{\Sigma^1}} A I_{[h]}, \bigoplus_{[k]\in\Sigma^1/\sim_{\Sigma^1}} A I_{[k]}\big]=\big[I_{[g]},  A I_{[g]}, A I_{[g]}\big]=0.
$$ 
Taking into account $A\LL=\LL,$ we get $I_{[g]}\subset Z(\LL) = 0,$ which is  a contradiction. Now, by $A=\bigoplus_{[\lam]\in\Lam^1/\approx_{\Lam^1}}\aa_{[\lam]},$ there exists $[\lam]\in\Lam^1/\approx_{\Lam^1}$ such that $\aa_{[\lam]} I_{[g]} \neq 0.$
 
Second, we will prove the uniqueness. Suppose that there exist $[\lam],~[\mu]\in\Lam^1/\approx_{\Lam^1}$ such that $\aa_{[\lam]} I_{[g]} \neq 0$ and $\aa_{[\mu]} I_{[g]} \neq 0$ for any $[g]\in\Sigma^1/\sim_{\Sigma^1}.$ From here, we can take $\lam_1\in[\lam],~\mu_2\in[\mu]$ and $g_1, g_2\in[g]$ such that $\aa_{\lam_1}I_{g_1}\neq 0$ and $\aa_{\mu_2}I_{g_2}\neq 0.$ Since $g_1, g_2\in[g],$ we can consider the $\Sigma^1$-connection $\{g_1, h_2, h_3, ..., h_{2n+1}\}\subset\Sigma\cup\Lam\cup\{1\}$ from $g_1$ to $g_2.$ We continue the proof in four cases;

{\bf Case 1.}  If $\lam_1 g_1\neq 1$ and $\mu_2 g_2\neq 1,$ then $\lam_1 g_1,~\mu_2 g_2\in\Lam^1.$ We can consider a $\Lam^1$-connection 
$$
\{\lam_1, g_1, {\lam_1}^{-1}, h_2, ..., h_{2n+1}, \mu_2, {g_2}^{-1}\}\subset\Sigma\cup\Lam,
$$ 
from $\lam_1$ to $\mu_2$  in the case $g_1h_2...h_{2n+1}=g_2,$ and  in  case $g_1h_2...h_{2n+1}={g_2}^{-1}$  is 
$$
\{\lam_1, g_1, {\lam_1}^{-1}, h_2, ..., h_{2n+1}, {\mu_2}^{-1}, g_2\}\subset\Sigma\cup\Lam.
$$
Then $\lam_1\approx_{\Lam^1}\mu_2$ and so $[\lam]=[\mu].$

{\bf Case 2.} If $\lam_1 g_1= 1$ and $\mu_2 g_2\neq 0,$ then $\lam_1={g_1}^{-1}$ and $\mu_2 g_2\in\Sigma^1.$ So we have a $\Lam^1$-connection 
$$
\{{g_1}^{-1},  {h_2}^{-1}, ..., {h_{2n+1}}^{-1}, {\mu_2}^{-1}, g_2\}\subset\Sigma\cup\Lam,
$$
from $\lam_1$ to $\mu_2$ in the case $g_1h_2...h_n=g_2.$ In the case $g_1h_2...h_n={g_2}^{-1}$ the $\Lam^1$-connection is
$$
\{{g_1}^{-1},  {h_2}^{-1}, ..., {h_{2n+1}}^{-1}, \mu_2, {g_2}^{-1}\}\subset\Sigma\cup\Lam.
$$
Then $\lam_1\approx_{\Lam^1}\mu_2,$ and so $[\lam]=[\mu].$

{\bf Case 3.} If $\lam_1 g_1\neq 1$ and $\mu_2 g_2= 1,$ then by a similar argumen as the second case we get $[\lam]=[\mu].$

{\bf Case 4.} If $\lam_1 g_1= 1$ and $\mu_2 g_2= 1,$ then $\lam_1={g_1}^{-1}$ and $\mu_2= {g_2}^{-1}.$ So we have a $\Lam^1$-connection 
$$
\{{g_1}^{-1},  {h_2}^{-1}, ..., {h_{2n+1}}^{-1}\}\subset\Sigma\cup\Lam.
$$
from $\lam_1$ to $\mu_2,$ and so $[\lam]=[\mu].$

Therefore, we conclude that for any $[g]\in\Sigma^1/\sim_{\Sigma^1}$
there exists a unique $[\lam]\in\Lam^1/\approx_{\Lam^1}$ such that $ \aa_{[\lam]} I_{[g]} \neq 0.$ \qed\\

It could be remarked that  Proposition \ref{5.10} shows that $I_{[g]}$
is an $A_{[\lam]}$-module. Hence we can assert the following result.

\begin{thm}\label{5.19}  Let $(\LL, A)$ be a tight graded $3$-Lie-Rinehart  algebra. Then
$$
\LL=\bigoplus_{i\in I}\LL_i,~~A=\bigoplus_{j\in J} A_j,
$$
where any $\LL_i$ is a non-zero graded  ideal of $\LL$ satisfying $[\LL_i, \LL_k]=0,$ when $i\neq k,$ and any $ A_j$ is a non-zero graded ideal of $A$ such that $A_j A_l=0$ when $j\neq l.$ Moreover, both decompositions satisfy that for any $i\in I$ there exists a unique $j\in J$ such that 
$$
A_j \LL_i\neq 0. 
$$
Fortheremore, any $(\LL_i, A_j)$ is a graded $3$-Lie-Rinehart  algebra. 
\end{thm}

===========================

\section{the  simple components of graded $3$-Lie-Rinehart algebras} 
\setcounter{equation}{0}\

In this section we focus on the simplicity of  graded $3$-Lie-Rinehart  algebra
$(\LL, A)$ by centering our attention in those of maximal length. From now on, we will suppose $\Sigma^1$ is symmetric, that is, if $g\in\Sigma^1$ then $g^{-1}\in\Sigma^1$ and also that $\Lam^1$ is symmetric in the same sence.

Let us introduce the concepts of root-multiplicativity and
maximal length in the framework  of  graded $3$-Lie-Rinehart   algebra, in a similar way to the ones for  split
Lie-Rinehart algebra in \cite{Kh}.

\begin{DEF}\label{mal} A graded  $3$-Lie-Rinehart   algebra
$(\LL, A)$  is called  {\em $G$-multiplicative} if for any $g, h, k\in\Sigma^1$ and $\lam, \mu\in\Lam^1$ the following conditions hold 
\begin{itemize}
	\item[-] If  $ghk\in\Sigma^1$\ then $[\LL_g, \LL_h, \LL_h]\neq 0.$
	\item[-]  If  $\lam g\in\Sigma^1$ then $A_\lam\LL_g\neq 0.$
	\item[-]  If  $\lam\mu\in\Lam^1$ then $A_\lam A_\mu\neq 0.$
\end{itemize}
\end{DEF}

\begin{DEF} A graded  $3$-Lie-Rinehart  algebra
	$(\LL, A)$  is called of {\em maximal length} if  for
any $g\in\Sigma^1$ and  $\lam\in\Lam^1$  we have $\dim
\LL_{g}=1=\dim A_\lam.$ 
\end{DEF}

\begin{rem} If $(\LL, A)$ is a graded  $3$-Lie-Rinehart  algebra such that $\LL$ and $A$ are ge-simple algebras then $Z(\LL) =\{0\}= Ann(A)$  and $Ann_{\LL}(A) = \{0\}.$  Also as consequence of Theorem \ref{main1}-(2) and Theorem \ref{main1'}-(2) we get that all of the non-zero elements
in $\Sigma^1$ are connected, that all of the non-zero elements in $\Lam^1$ are also connected and that
\begin{eqnarray*}
\LL_1=\sum_{g\in\Sigma^1\cap\Lam^1} A_{g^{-1}}\LL_g
+\sum_{h, k\in\Sigma^1}[\LL_h, \LL_k, \LL_{(hk)^{-1}}],\\
A_1=\sum_{\mu\in\Lam^1} A_{\mu^{-1}}A_\mu+\sum_{h, k\in {\Lam^1\cap\Sigma^1}}\rho(\LL_h, \LL_k) A_{(hk)^{-1}}.
\end{eqnarray*}
From here, the conditions for $(\LL, A)$ of being tight 
together with the ones of having $\Sigma^1$ and $\Lam^1$ all of their elements connected, are necessary
conditions to get a characterization of the gr-simplicity of the algebras $\LL$ and $A.$ Actually, we are going to shwo that 
under the hypothesis of being $(\LL, A)$ of maximal length and $G$-multiplicative, these are
also sufficient conditions. 
\end{rem}

\begin{lem}\label{10}  Let $(\LL, A)$ be a tight  graded $3$-Lie-Rinehart  algebra of maximal length and $G$-multiplicative. If $I$ is a  graded ideal of $\LL$ such that $I\subset\LL_1,$ then  $I=\{0\}.$ 
\end{lem}
\noindent {\bf Proof.} Suppose there exists a non-zero graded ideal $I$  of $\LL$ such that $I\subset\LL_1.$  We are going to show that $I\subset Ann_{\LL}(A).$ 

If $g, h\in\Sigma^1$ with $gh\neq 1,$ we have 
 \begin{equation}\label{502}
[I, \LL_g, \LL_h]\subset\LL_{gh}\cap\LL_1=0.		
\end{equation}
Otherwise, if $h=g^{-1}$ and $[I, \LL_g, \LL_{g^{-1}}]\neq 0$ for some $g\in\Sigma^1,$ then there exist $x\in\LL_g, x'\in\LL_{g^{-1}}$ and $i\in I$ such that $[i, x, x']\neq 0.$ By the $G$-maltiplicativity (consider $1, g, 1\in\Sigma^1\cup\{1\}$) and  the maximal length of $(\LL, A)$ there exist $x_1\in\LL_1$ such that $0\neq[i, x, x_1]\in I$ a contradiction. Hence,
\begin{equation}\label{503}
[I, \LL_g, \LL_{g^{-1}}]=0,~~~~~\forall g\in\Sigma^1.		
\end{equation}
By Eqs. (\ref{502}) and (\ref{503}), we get $[I, \LL_g, \LL_h]=0$ for all  $g, h\in\Sigma^1.$ therefore,
\begin{equation}\label{504}
[I, \LL, \LL]=0.		
\end{equation}
Next, we show that $A I=0.$ Note that
\begin{equation}\label{510}
A I=\bigg(A_1\oplus\big(\bigoplus_{\lam\in \Lambda^1}A_\lam\big)\bigg)I\subset A_1 I+ \bigoplus_{\lam\in \Lambda^1}A_\lam I.		
\end{equation}
For the second summand in (\ref{510}), since $(\LL, A)$ is a tight  graded $3$-Lie-Rinehart  algebra,  we have
\begin{equation}\label{511}
 \bigoplus_{\lam\in \Lambda^1}A_\lam I\subset  \bigoplus_{\lam\in \Sigma^1}\LL_\lam\cap\LL_1=0.		
\end{equation}
Now, consider the first summand in (\ref{510}), since $$A_1=\sum_{\mu\in\Lam^1} A_{\mu^{-1}}A_\mu+\sum_{h, k\in {\Lam^1\cap\Sigma^1}}\rho(\LL_h, \LL_k) A_{(hk)^{-1}},$$ we have
\begin{equation}\label{512}
 A_1 I\subset\sum_{\mu\in\Lam^1}\big( A_{\mu^{-1}}A_\mu\big)I+\sum_{h, k\in {\Lam^1\cap\Sigma^1}}\rho(\LL_h, \LL_k) \big(A_{(hk)^{-1}}\big)I.		
\end{equation}
For the first item in (\ref{512}), by the fact that $I$ is an $A$-module and Eq. (\ref{511}), we get 
\begin{equation}\label{513}
 ( A_{\mu^{-1}}A_\mu)I=A_{\mu^{-1}}(A_\mu I)=0.		
\end{equation}
Consider the second item in (\ref{512}), by Eq.  (\ref{4}),  the fact that $I$ is an $A$-module and Eq. (\ref{504}), we have
\begin{eqnarray}\label{514}
\nonumber\rho(\LL_h, \LL_k) (A_{(hk)^{-1}})I&\subset&[\LL_h, \LL_k, A_{(hk)^{-1}}I]+ A_{(hk)^{-1}}[\LL_h, \LL_k, I]\\
&\subset&[\LL_h, \LL_k, I]+ A_{(hk)^{-1}}[\LL_h, \LL_k, I]\\
\nonumber&=&0.
\end{eqnarray}
 Eqs.  (\ref{513}) and  (\ref{514}), give us
\begin{equation}\label{515}
 A_1 I=0.		
\end{equation}
Now, Eqs.  (\ref{511}) and  (\ref{515}), implies that $A I=0,$ taking into account Eq. (\ref{504}) we obtain   $I\subset Ann_{\LL}(A)=0.$\qed

\begin{pro}\label{6.4} Let $(\LL, A)$ be a tight  graded $3$-Lie-Rinehart  algebra of maximal length and $G$-multiplicative. If all the elements in $\Sigma^1$ are $\Sigma^1$-connected,  then either $\LL$ is gr-simple or $\LL=I\oplus I'$ where  $I$ and $I'$ are graded simple ideals of $\LL.$ 
\end{pro}	
\noindent {\bf Proof.}  Consider $I$ a nonzero graded ideal of $\LL.$ By Lemma \ref{10}, we have   $I\nsubseteq\LL_1$ and  the maximal length of $\LL$ gives us

\begin{equation*}
I=(I\cap\LL_1)\oplus\big(\bigoplus_{g\in\Sigma^1}(I\cap\LL_g)\big),
\end{equation*}
with $(I\cap\LL_g)\neq 0$ for some $g\in\Sigma^1.$ Denote by $I_g~:=I\cap\LL_g$ and by $\Sigma^1_I :=\{g\in\Sigma^1 : I_g\neq 0\}=\{g\in\Sigma^1 : I\subset\LL_g\}.$  Then we can rewrite 
$$
I=(I\cap\LL_1)\oplus\big(\bigoplus_{g\in\Sigma^1_I}I_g\big),
$$ 
with $\Sigma^1_I\neq\emptyset.$  Let us distinguish two cases.

{\bf Case 1.} Suppose there exists $g_0\in\Sigma^1_I$ such that ${g_0}^{-1}\in\Sigma^1_I.$  Then $0\neq I_{g_0}\subset I$ and  by the maximal length of $(\LL, A)$  we have
\begin{equation}\label{520}
0\neq\LL_{g_0}\subset I.
\end{equation}

Now, let us take some $h\in\Sigma^1$ satisfying $h\notin\{g_0, {g_0}^{-1}\}.$ By  the assumption, $g_0$ is $\Sigma^1$-connected to $h,$ then  we have a $\Sigma^1$-connection  $\{g_1, g_2,...,
g_{2n+1}\}\subset\Sigma^1\cup\Lam^1\cup\{1\}.$ 

Consider $g_1, g_2, g_3\in\Sigma^1\cup\Lam^1\cup\{1\}$ and if $g_1g_2g_3\in\Sigma$ (respectively, $g_1g_2\in\Lam$),  since $g_0=g_1\in\Sigma^1_I$  we hav $\LL_{g_1}\neq 0.$ From here, the $G$-multiplicativity and maximal length of $\LL$ allow us to get
$$
0\neq[\LL_{g_1}, \LL_{g_2}, \LL_{g_3}]=\LL_{g_1g_2g_3}~~\hbox{~~(respectively}, 0\neq A_{g_1}\LL_{g_2}=\LL_{g_1g_2}).
$$
Since $0\neq\LL_{g_1}\subset I$ as consequence of Eq. (\ref{520}),  we have 
$$
0\neq\LL_{g_1g_2g_3}\subset I ~~\hbox{~~(respectively},  \LL_{g_1g_2}\subset I).
$$
We can follow this process with the connection $\{g_1, g_2,...,g_{2n+1}\}$  and obtain that
$$
0\neq\LL_{g_1g_2g_3 ...g_{2n+1}}\subset I.
$$
Thus we have shown that 
\begin{equation}\label{521}
\hbox{for ~~any~~} h\in\Sigma^1,~ \hbox{we~~ have~~ that~}~ 0\neq \LL_k\subset I\hbox{~~for~~some~~}~~k\in\{h, h^{-1}\}.
\end{equation}	
Since $g_0^{-1}\in\Sigma^1_I,$  we  have $\{g^{-1}_1, g^{-1}_2, ..., g^{-1}_{2n+1}\}$ is a $\Sigma^1$-connection from $g_0^{-1}$ to $h$ satisfying
$$
g^{-1}_1 g^{-1}_2g_3^{-1} ...g^{-1}_{2n+1}=k^{-1}.
$$
By arguing as above we get,
\begin{equation}\label{522}
0\neq \LL_{k^{-1}}\subset I,	
\end{equation}	
and so $\Sigma^1_I=\Sigma^1.$ The fact 
$
\LL_1=\sum_{g\in\Sigma^1\cap\Lam^1} A_{g^{-1}}\LL_g
+\sum_{h, k\in\Sigma^1}[\LL_h, \LL_k, \LL_{(hk)^{-1}}],
$
implies that
\begin{equation}\label{523}
\LL_1\subset I.
\end{equation}
From Eqs. (\ref{520})-(\ref{523}),  we obtain $\LL\subset I,$  and so $\LL$ is gr-simple.

{\bf Case 2.} In the second case, suppose that for any $g_0\in\Sigma^1_I$ we have that $g_0^{-1}\notin\Sigma^1_I.$ Observe that by
arguing as in the  case 1,  we can write
\begin{equation}\label{530}
\Sigma^1=\Sigma^1_I\cup\Sigma_I^c,
\end{equation}
where $\Sigma_I^c=\{g^{-1}~:~ g\in\Sigma^1_I\}.$ Denote by
$$
I'~:=\sum_{g\in\Lam^1,~ g^{-1}\in \Sigma_I^c}A_g\LL_{g{-1}}\oplus\big(\bigoplus_{g'\in\Sigma_I^c}\LL_{g'}\big).
$$
We are giong to show that $I'$ is a graded ideal  of $3$-Lie-Rinehart algebra  $(\LL, A).$ By construction $I'$ is $G$-graded.  First, we will show that $I'$ is a $3$-Lie ideal of $\LL.$ Taking into account Eq. (\ref{00}), we have
\begin{eqnarray}\label{531}
\nonumber[\LL, \LL, I']&=&\bigg[\LL_1\oplus(\bigoplus_{h\in\Sigma^1}\LL_h), \LL_1\oplus(\bigoplus_{k\in\Sigma^1}\LL_k), \sum_{g\in\Lam^1,~ g^{-1}\in \Sigma_I^c}A_g\LL_{g{-1}}\oplus(\bigoplus_{g'\in\Sigma_I^c}\LL_{g'})\bigg]\\
&\subset&\big[\LL_1, \LL_1, \sum_{g\in\Lam^1,~ g^{-1}\in \Sigma_I^c}A_g\LL_{g{-1}}\big]+\big[\LL_1, \LL_1, \bigoplus_{g'\in\Sigma_I^c}\LL_{g'}\big]\\
\nonumber&+&\big[\LL_1, \bigoplus_{h\in\Sigma^1}\LL_h, \sum_{g\in\Lam^1,~ g^{-1}\in \Sigma_I^c}A_g\LL_{g{-1}}\big]+\big[\LL_1, \bigoplus_{k\in\Sigma^1}\LL_k,  \bigoplus_{g'\in\Sigma_I^c}\LL_{g'}\big]\\
\nonumber&+&\big[\bigoplus_{h\in\Sigma^1}\LL_h, \LL_1, \sum_{g\in\Lam^1,~ g^{-1}\in \Sigma_I^c}A_g\LL_{g{-1}}\big]+\big[\bigoplus_{h\in\Sigma^1}\LL_h, \LL_1, \bigoplus_{g'\in\Sigma_I^c}\LL_{g'}\big]\\
\nonumber&+&\big[\bigoplus_{h\in\Sigma^1}\LL_h, \bigoplus_{k\in\Sigma^1}\LL_k, \sum_{\substack{g\in\Lam^1\\ g^{-1}\in \Sigma_I^c}}A_g\LL_{g{-1}}\big]+\big[\bigoplus_{h\in\Sigma^1}\LL_h, \bigoplus_{k\in\Sigma^1}\LL_k, \bigoplus_{g'\in\Sigma_I^c}\LL_{g'}\big].
\end{eqnarray}
For the first summand in (\ref{531}), if there exist $g\in\Lam^1$ and $g^{-1}\in \Sigma_I^c$ such that $[\LL_1, \LL_1, A_g\LL_{g{-1}}]\neq 0,$ by Eq.(\ref{4}) we have
\begin{eqnarray*}
[\LL_1, \LL_1, A_g\LL_{g{-1}}]&=&A_g[\LL_1, \LL_1, \LL_{g{-1}}]+\rho(\LL_1, \LL_1) A_g\LL_{g{-1}}\\
&\subset& A_g\LL_{g{-1}}\subset I'.
\end{eqnarray*}
Therefore,
\begin{equation}\label{532}
[\LL_1, \LL_1, \sum_{g\in\Lam^1,~ g^{-1}\in \Sigma_I^c}A_g\LL_{g{-1}}]\subset I'. 
\end{equation}
For the second summand in  (\ref{531}), it is clear that 
\begin{equation}\label{533}
[\LL_1, \LL_1, \bigoplus_{g'\in\Sigma_I^c}\LL_{g'}]\subset I'.
\end{equation}
Connsider the third summand in  (\ref{531}), if $[\LL_1, \LL_h, A_g\LL_{g{-1}}]\neq 0$ for some $g\in\Lam^1,  g^{-1}\in \Sigma_I^c$ and $h\in\Sigma^1.$ Then in case $h=g^{-1}$ clearly $[\LL_1, \LL_h, A_g\LL_{g{-1}}]\subset\LL_h\subset I',$ and in case $h=g,$  the maximal length of $\LL$ and the fact $I$ is a graded ideal  give us
$$
\LL_h=[\LL_1, \LL_h, A_g\LL_{g^{-1}}]\subset I\cap I'=\{0\},
$$ 
which is a contradiction with $g\in\Sigma^1_I.$ Now, if $h\notin\{g, g^{-1}\},$ we then have
$$
0\neq [\LL_1, \LL_h, A_g\LL_{g^{-1}}]\subset A_g[\LL_1, \LL_h, \LL_{g^{-1}}]+\rho(\LL_1, \LL_h)A_g\LL_{g^{-1}}.
$$
By the maximal length of $\LL,$   
$$
\hbox{ either} ~~0\neq  A_g[\LL_1, \LL_h, \LL_{g^{-1}}]=\LL_h~~\hbox{ or}~~~ 0\neq\rho(\LL_1, \LL_h)A_g\LL_{g^{-1}}=\LL_h.
$$
In both cases, by $G$-multiplicativity, we have that $\LL_{h^{-1}}\subset I$ and therefore $h^{-1}\in\Sigma^1_I,$ this implies that $h\in\Sigma^c_I$ and then $\LL_h\subset I'.$ Hence,
\begin{equation}\label{534}
[\LL_1, \bigoplus_{h\in\Sigma^1}\LL_h, \sum_{g\in\Lam^1,~ g^{-1}\in \Sigma_I^c}A_g\LL_{g{-1}}]\subset I'. 
\end{equation}
By the skew symmetry,
\begin{equation}\label{535}
[\bigoplus_{h\in\Sigma^1}\LL_h, \LL_1, \sum_{g\in\Lam^1,~ g^{-1}\in \Sigma_I^c}A_g\LL_{g{-1}}]\subset I'. 
\end{equation}
A similar argument as above for the seventh summand in (\ref{531}), one can show that
\begin{equation}\label{536}
\big[\bigoplus_{h\in\Sigma^1}\LL_h, \bigoplus_{k\in\Sigma^1}\LL_k, \sum_{\substack{g\in\Lam^1\\ g^{-1}\in \Sigma_I^c}}A_g\LL_{g{-1}}\big]\subset I'.
 \end{equation}
Next, consider the fourt summand in  (\ref{531}), suppose  there exist $k\in\Sigma^1$ and $g'\in\Sigma_I^c$ such that $[\LL_1, \LL_k, \LL_{g'}]\neq 0.$ In case $k=g'^{-1},$ we have $0\neq[\LL_1, \LL_k, \LL_{g'}]\subset I.$ Now, since $I$  is a graded ideal and  $\LL$ is $G$-multiplicative,  we have $$\LL_{g'}=\big[[\LL_1, \LL_k, \LL_{g'}], \LL_1, \LL_{g'}\big]\subset I,$$ and so $g'\in\Sigma^1_I$ a contradiction with   $g'\in\Sigma_I^c.$ In case $k\neq g'^{-1},$  the $G$-multiplicativity gives us $\LL_{k^{-1}{{g'}^{-1}}}=[\LL_1, \LL_{k^{-1}}, \LL_{{g'}^{-1}}]\subset I.$ From here $k^{-1}{{g'}^{-1}}\in\Sigma^1_I$ and so $kg'\in\Sigma_I^c.$ Thus, we get $[\LL_1, \LL_k, \LL_{g'}]=\LL_{kg'}\subset I'.$ Therefor,
\begin{equation}\label{537}
\big[\LL_1, \bigoplus_{k\in\Sigma^1}\LL_k,  \bigoplus_{g'\in\Sigma_I^c}\LL_{g'}\big]\subset I'. 
\end{equation}
By the skew symmetry,
\begin{equation}\label{538}
\big[\bigoplus_{h\in\Sigma^1}\LL_h, \LL_1, \bigoplus_{g'\in\Sigma_I^c}\LL_{g'}\big]\subset I'. 
\end{equation}
Finally, for the last summand in (\ref{531}). Suppose $0\neq[\LL_h, \LL_k, \LL_{g'}]$ for some $h, k\in\Sigma^1$ and $g'\in\Sigma_I^c.$ If $hk=1,$ clearly $[\LL_h, \LL_k, \LL_{g'}]=\LL_{g'}\subset I'.$ Now, if $g'\neq h^{-1}$ and $g'\neq k^{-1},$   the $G$-multiplicativity and maximal length of $\LL$ allow us to get $\LL_{g'}=[\LL_h, \LL_k, \LL_{g'}]\subset I,$ a contradiction. In case $g'\neq h^{-1},$ we have $(hk)^{-1}$ and so $\LL_{g'}=[\LL_h, \LL_k, \LL_{g'}]\subset I.$ 
Therefore,
\begin{equation}\label{539}
\big[\bigoplus_{h\in\Sigma^1}\LL_h, \LL_1, \bigoplus_{g'\in\Sigma_I^c}\LL_{g'}\big]\subset I'. 
\end{equation}
From Eqs. (\ref{532})-(\ref{539}), we conclud that $I'$ is a $3$-Lie ideal of $\LL.$

Second, we will check $A I'\subset I'.$   We have

\begin{eqnarray}\label{370}
\nonumber A I'&=&\bigg(A_1\oplus\big(\bigoplus_{\lam\in \Lambda^1}A_\lam\big)\bigg)\bigg(\sum_{g\in\Lam^1,~ g^{-1}\in \Sigma_I^c}A_g\LL_{g{-1}}\oplus\big(\bigoplus_{g'\in\Sigma_I^c}\LL_{g'}\big)\bigg)\\
&\subset&I'+\big(\bigoplus_{\lam\in \Lambda^1}A_\lam\big)\big(\sum_{g\in\Lam^1,~ g^{-1}\in \Sigma_I^c}A_g\LL_{g{-1}}\big)+\big(\bigoplus_{\lam\in \Lambda^1}A_\lam\big)\big(\bigoplus_{g'\in\Sigma_I^c}\LL_{g'}\big)
\end{eqnarray}
Consider the third summand in (\ref{370}) and suppose that $A_\lam \LL_{g'}\neq 0$ for some $\lam\in\Lam^1,~g'\in\Sigma^c_I.$ If $\lam g'\in\Sigma^1_I,$ so $\lam^{-1}{g'^{-1}}\in\Sigma^1$ then by the $G$-multiplicativity of $\LL$ we get $A_{\lam^{-1}}\LL_{g'^{-1}}\neq 0.$  Now by the maximal length of $\LL$ and the fact ${g'^{-1}}\in\Sigma^1_I,$ we conclud that $A_{\lam^{-1}}\LL_{g'^{-1}}=\LL_{\lam^{-1}g'^{-1}}\subset I.$ Therefore ${\lam^{-1}g'^{-1}}=(\lam g')^{-1}\in\Sigma^1_I$ which is a contradiction. Hence $\lam g'\in\Sigma^c_I,$ and so $A_\lam \LL_{g'}\subset I'.$ Therefore,
\begin{equation}\label{363}
\big(\bigoplus_{\lam\in \Lambda^1}A_\lam\big)\big(\bigoplus_{g'\in\Sigma_I^c}\LL_{g'}\big)\subset I'.
\end{equation}
We can argue as above with the second summand in (\ref{370}),  so as to conclude that
\begin{equation}\label{375.}
\big(\bigoplus_{\lam\in \Lambda^1}A_\lam\big)\big(\sum_{g\in\Lam^1,~ g^{-1}\in \Sigma_I^c}A_g\LL_{g{-1}}\big)\subset I'.
\end{equation}
From Eqs. (\ref{363}) and (\ref{375.}) we get $A I'\subset I'.$

Finally, let us check $\rho\big(I', I'\big)(A)\LL\subset I'.$ In fact by Eq. (\ref{4}) we have
$$
\rho\big(I', I'\big)(A)\LL\subset [I', I', A\LL]+A[I', I', \LL]
$$
Tanks to $I'$ is a $3$-Lie  ideal we get the result. 

Summarizing a discussion of above, we conclude that $I'$ is a graded ideal of the graded $3$-Lie-Rinehart  algebra $(\LL, A).$

Next, by Eq. (\ref{530}) we get $\sum_{h, k\in\Sigma^1}[\LL_h, \LL_k, \LL_{(hk)^{-1}}]=0,$ so by hypothesis we  must have 
$$
\LL_1=\sum_{g\in\Sigma^1\cap\Lam^1} A_{g^{-1}}\LL_g=\sum_{g\in\Sigma^1_I,~g^{-1}\in\Lam^1}A_{g^{-1}}\LL_g\oplus\sum_{g^{-1}\in\Sigma_I^c,~g\in\Lam^1}A_g\LL_{g^{-1}}.
$$
For direct character, take
$$
0\neq x\in\sum_{g\in\Sigma^1_I,~g^{-1}\in\Lam^1}A_{g^{-1}}\LL_g\cap\sum_{g^{-1}\in\Sigma_I^c,~g\in\Lam^1}A_g\LL_{g^{-1}}.
$$
Taking into account $Z_\rho(\LL) = \{0\}$ and $\LL$ is graded, there exist $0\neq y\in\LL_h,~ 0\neq z\in\LL_k$ for some $h, k\in\Sigma^1$ such that $[x, y, z]\neq 0,$ being then $\LL_h\in I\cap I'=\{0\},$ a contradiction. Hence  the sum is direct. Taking into account the above observation and Eq. (\ref{530}) we
have
$$
\LL=I\oplus I'.
$$
Note that, one can proceed with $I$ and $I'$ as we did for $\LL$ in the first case of the proof to
conclude that $I$ and $I'$ are graded simple ideals of $\LL,$  which completes the proof of the proposition.\qed\\

In a similar way to Proposition \ref{6.4},  one can prove the next result;

\begin{pro}\label{6.4.} Let $(\LL, A)$ be a tight  graded  $3$-Lie-Rinehart   algebra of maximal length and $G$-multiplicative. If all the elements in $\Lam^1$ are $\Lam^1$-connected,  then either $A$ is gr-simple or $A=J\oplus J'$ where  $J$ and $J'$ are graded simple ideals of $A.$\qed\\ 
\end{pro}	

Now, we are ready to state our main result;

\begin{thm}\label{final}  Let $(\LL, A)$ be a tight graded $3$-Lie-Rinehart  algebra of maximal length, $G$-multiplicative, with symmetric $G-$supports $\Sigma^1$ and $\Lam^1$ in such a way that  $\Sigma^1$  and $\Lam^1$ have  all their elements $\Sigma^1$-connected and  $\Lam^1$-connected, respectively.Then
$$
\LL=\bigoplus_{i\in I}\LL_i,~~\hbox{~and~~}~~ A=\bigoplus_{j\in J} A_j,
$$
where any $\LL_i$ is a graded simple  ideal of $\LL$  having all of its elements in $G$-support  $\Sigma^1$-connected and such
that $[\LL_{i_1}, \LL_{i_2}, \LL_{i_3}]=0$ for any $i_1, i_2, i_3\in I$  different from each other,  and any $A_j$ is a graded simple ideal of $A$  satisfying $A_j A_l=0$ for any $l\in J$ such that $j\neq l.$  Moreover, both decompositions satisfy that for any $r\in I$ there exists a unique $\bar{r}\in J$ such that 
$$
A_{\bar{r}} \LL_r\neq 0. 
$$
Fortheremore, any $(\LL_r, A_{\bar{r}}, \rho|_{\LL_i\times\LL_i})$ is a graded $3$-Lie-Rinehart  algebra. 
\end{thm}
\noindent {\bf Proof.} By Theorem \ref{5.19}  we can  write
$$
\LL=\bigoplus_{[g]\in\Sigma^1/\sim_{\Sigma^1}}I_{[g]},
$$
with any $I_{[g]}$ a graded ideal of $\LL,$ being each $I_{[g]}$ a graded $3$-Lie-Rinehart  algebra having as $G$-support $[g].$  Also we can write $A$ as the direct sum of the graded ideals
$$
 A=\bigoplus_{[\lam]\in\Lam^1/\approx_{\Lam^1}}\aa_{[\lam]},
$$
in such a way that any $\aa_{[\lam]}$ has as $G$-support $[\lam],$ for any $[g]\in\Sigma^1/\sim_{\Sigma^1}$ there exists a unique $[\lam]\in\Lam^1/\approx_{\Lam^1}$ such that $ \aa_{[\lam]} I_{[g]} \neq 0$ and being $( I_{[g]}, \aa_{[\lam]})$ a graded $3$-Lie-Rinehart  algebra.

Now, by applying Proposition \ref{6.4} and Proposition \ref{6.4.} to each $( I_{[g]}, \aa_{[\lam]}),$  in a similar manner to observe that the $\Sigma^1$-multiplicativity of  $( I_{[g]}, \aa_{[\lam]}),$ that is, $( I_{[g]}, \aa_{[\lam]})$ is $\Sigma^1$-multiplicative as consequence of the $\Sigma^1$-multiplicativity of $(\LL, A).$  Clearly  $( I_{[g]}, \aa_{[\lam]})$ is of maximal length. We also have $( I_{[g]}, \aa_{[\lam]})$ is tight,  as consequence of tightness of $(\LL, A)$ (see Proposition \ref{6.4} and Proposition \ref{6.4.}). 

Next,  we can apply Proposition \ref{6.4} and Proposition \ref{6.4.} to each $( I_{[g]}, \aa_{[\lam]})$ so as to conclude that any $ I_{[g]}$
is either graded simple or the direct sum of graded simple
ideals $ I_{[g]}=P\oplus Q,$ and that any $\aa_{[\lam]}$
is either graded simple or the direct sum of graded simple ideals $\aa_{[\lam]}=R\oplus S.$ From here, it is clear that by writing $\LL_i=P\oplus Q$ and $\aa_j=R\oplus S$ if $\LL_i$ or $\aa_j$ are not  graded simple. Then Theorem \ref{5.19} allows as to assert that the resulting decomposition satisfies the assertions of the theorem. \qed\\



\end{document}